\documentclass[final,hidelinks,onefignum,onetabnum]{siamart251216}

\usepackage{amsfonts,amssymb,amsmath,booktabs,multirow}
\usepackage{xcolor}         
\usepackage[sort,nocompress]{cite}

 \usepackage{placeins}
\usepackage{dsfont}
\usepackage{graphicx}
\usepackage{epstopdf}
\usepackage{algorithmic}

\ifpdf
  \DeclareGraphicsExtensions{.eps,.pdf,.png,.jpg}
\else
  \DeclareGraphicsExtensions{.eps}
\fi

\usepackage{amsmath}
\usepackage{graphicx}
\usepackage{listings}
\usepackage{subfigure}
\usepackage{subeqnarray}
\usepackage{amsfonts}
\usepackage{amssymb}
\usepackage{amsbsy}
\usepackage{amsmath}
\usepackage[mathscr]{eucal}
\usepackage{color}
\usepackage{lipsum}
\usepackage{amsfonts}
\usepackage{graphicx}
\usepackage{epstopdf}
\usepackage{algorithmic}
\usepackage{algorithm}
\usepackage{amsopn}
\usepackage{makecell}
\usepackage{multirow}

\usepackage[font=scriptsize, labelfont=normalfont, textfont=normalfont, labelsep=colon, justification=justified, singlelinecheck=true]{caption}

\newcommand{\customref}[2]{\hyperref[#2]{#1}}

\newtheorem{remark}{Remark}[section]
\newtheorem{assumption}{Assumption}
\crefname{assumption}{Assumption}{Assumptions}
\makeatletter
\@ifpackageloaded{lineno}{%
  \@ifpackagelater{lineno}{2023/05/18}{%
    \allowdisplaybreaks[4]%
  }{%
    \PackageWarningNoLine{siam_template}{%
      lineno older than v5.2 detected; disabling
      \string\allowdisplaybreaks\space to preserve pagination%
    }%
  }%
}{%
  \allowdisplaybreaks[4]%
}
\makeatother

\usepackage{mathtools}

\usepackage[normalem]{ulem} 
\usepackage{dsfont} 
\usepackage{xcolor} 
\usepackage{enumitem}

\newcommand{\lb}{\lambda}

\newcommand{\cB}{{\mathcal{B}}}
\newcommand{\cC}{{\mathcal{C}}}

\newcommand{\cP}{{\mathcal{P}}}

\newcommand{\cX}{{\mathcal{X}}}

\newcommand{\RR}{\mathbb{R}} 

\newcommand{\dist}{\mathrm{dist}}    

\newcommand{\st}{\mbox{ s.t. }}

\DeclareMathOperator*{\argmin}{arg\,min} 

\newcommand{\bc}{\begin{center}}
\newcommand{\ec}{\end{center}}

\newcommand{\bdm}{\begin{displaymath}}
\newcommand{\edm}{\end{displaymath}}

\newcommand{\beq}{\begin{equation}}
\newcommand{\eeq}{\end{equation}}

\newcommand{\bfl}{\begin{flushleft}}
\newcommand{\efl}{\end{flushleft}}

\newcommand{\bt}{\begin{tabbing}}
\newcommand{\et}{\end{tabbing}}

\newcommand{\beqn}{\begin{eqnarray}}
\newcommand{\eeqn}{\end{eqnarray}}

\newcommand{\beqs}{\begin{align*}} 
\newcommand{\eeqs}{\end{align*}}  

\headers{First-order smoothed exact-penalty methods}{ }

\title{First-order smoothed exact-penalty methods for equality-constrained regularized nonconvex optimization
}

\author{Shiji Zuo\thanks{Department of Mathematical Sciences, Rensselaer Polytechnic Institute, Troy, NY (\email{zuos@rpi.edu, xuy21@rpi.edu}).}
\and Qihang Lin\thanks{Tippie College of Business, University of Iowa, Iowa City, IA (\email{qihang-lin@uiowa.edu}).}
\and Wei Liu\thanks{IMAI, Wuhan University; Department of Applied Mathematics, Hong Kong Polytechnic University, China (\email{lwdsdqqb@gmail.com}).}\and Yangyang Xu\footnotemark[1]} 

\usepackage{amsopn}

\ifpdf
\hypersetup{
  pdftitle={Smooth penalty method},
  pdfauthor={S. Zuo}
}
\fi

\begin{document}

\maketitle

\begin{abstract}  We consider regularized nonconvex optimization over a closed convex set $\mathcal{X} \subseteq \RR^d$, which involves a convex nonsmooth regularizer $r$ in the objective and $n$ smooth equality constraints. To find a near-KKT point of the problem, we propose a prox-linear first-order method, named as SEPi-Prox, based on a smoothed exact-penalty formulation. At each outer iteration of SEPi-Prox, a strongly convex prox-linear subproblem is formed and solved inexactly to a sufficient accuracy. 
Under a local non-singularity condition that is  
weaker than some other commonly used regularity conditions or constraint qualifications, SEPi-Prox finds an $\epsilon$-KKT point within $O(\epsilon^{-2})$ first-order oracle calls, matching the lower bound for unconstrained smooth nonconvex optimization. We employ 
Nesterov's accelerated proximal-gradient method (APG) and a dual cutting-plane method (DCP) for the \mbox{prox-linear} subproblems. Let $\mathcal{C}_{\mathrm{FO}}$ and $\mathcal{C}_r$ denote the numbers of arithmetic 
operations of evaluating one first-order oracle and one proximal mapping of $r$ over~$\mathcal{X}$, respectively. 
Then \mbox{SEPi-Prox} requires \mbox{$\widetilde{O}(\mathcal{C}_{\mathrm{FO}}\epsilon^{-2} + (dn + \mathcal{C}_r)\epsilon^{-5/2})$} and \mbox{$\widetilde{O}((\mathcal{C}_{\mathrm{FO}} + dn^2 + n\mathcal{C}_r + n^{3.373})\epsilon^{-2})$} total arithmetic operations to produce an $\epsilon$-KKT point, respectively by using APG and DCP as the subsolver. 
The latter result is preferable when $n$ is small and nearly matches the optimal dependence on $\epsilon$. 
Numerical experiments on instances of the multi-class Neyman--Pearson classification demonstrate the efficiency of \mbox{SEPi-Prox with both subsolvers.} 
\end{abstract}

 \begin{keywords}
 penalty method, nonconvex constrained optimization, composite optimization, first-order method, arithmetic complexity, cutting-plane method
 \end{keywords}

\begin{MSCcodes}
90C26, 90C30, 65K05, 49M37
\end{MSCcodes}

\section{Introduction}

	 In this paper, we focus on designing first-order methods (FOMs) for solving the \textit{nonconvex constrained problem}:  
 \begin{align}\label{eq:model}
			\min_{ x\in \mathcal{X} }   f(x) + r(x)
			, \text{ s.t. } g(x) = 0, 
			\tag{P} 
	\end{align}
	where $f: \RR^d \to \RR$ and $g: \RR^d \to \RR^n$ are smooth and possibly nonconvex mappings, $\mathcal{X}$ is a closed convex set, and $r: \RR^d \to \RR$ is a nonnegative and convex nonsmooth regularizer.  
   For simplicity, we formulate  \cref{eq:model} with equality constraints. Inequality constraints $q(x)\leq 0$ reduce to $q(x)+s=0$ by introducing slack variables $s\ge 0$. 
 Problem~\cref{eq:model} arises in various applications 
 such as Neyman--Pearson classification~\cite{rigollet2011neyman}, fairness-constrained classification~\cite{pmlr-v54-zafar17a}, phase retrieval problems~\cite{waldspurger2015phase}, and constrained deep learning~\cite{liu2023inexact}.  
 Throughout this paper, we impose the following assumption. 
	\begin{assumption}
		\label{ass:problemsetup}  
        The set $\mathcal X$ is closed and convex, and $r(\cdot)$ is convex  
        and nonnegative on $\mathcal X$.  
        There exist constants $l_r,l_g \ge 0$ such that  \mbox{$\max_{\xi \in \partial r(x)} \|\xi\|\le l_r$} and $\|J_g(x)\|\le l_g$ for all $x  \in \mathcal X$,
        where $J_g(x)\in\RR^{n\times d}$ denotes the Jacobian of $g$. 
        In addition, $f(x) \ge f^* > -\infty, \forall\, x\in \RR^d$, and there exist constants $L_f,L_g> 0$ such that 
        \begin{align*}
            \|\nabla f(x)-\nabla f(y)\|\le L_f\|x-y\|,\quad \|J_g(x)-J_g(y)\|\le L_g\|x-y\|,\quad \forall x, y \in \RR^d.
        \end{align*}
	\end{assumption}

	Because \eqref{eq:model} is nonconvex, finding a global minimizer is generally NP-hard, even for 
    quadratic programs~\cite{murty1985some}.  
    Hence, we aim to compute an $\epsilon$-KKT point of~\eqref{eq:model}, as defined in Definition~\ref{def:nearkkt}, for any target accuracy $\epsilon > 0$. Moreover,  
    the first-order oracle $(f,\nabla f, g, J_g)$  
    can be evaluated at any query point $x$ by at most $\mathcal{C}_{\mathrm{FO}}$ arithmetic operations. For any $c>0$, the proximal mapping of $c\cdot  r+\delta_{\mathcal{X}}$ requires at most $\mathcal{C}_{r}$ arithmetic operations, where $\delta_{\mathcal{X}}$ denotes the indicator function of~$\mathcal{X}$.  
   We measure oracle and arithmetic complexities by the first-order oracle calls and total arithmetic operations required to reach an $\epsilon$-KKT point. 

 FOMs have been extensively studied for functionally constrained convex optimization \cite{lin2018level-SIOPT,yu2017simple,xu2022first} and for minimizing nonconvex objectives subject to convex functional constraints \cite{huang2025inexact,zhang2020proximal}.  
For problem \eqref{eq:model}, the augmented Lagrangian method (ALM) in~\cite{li2021rate} and the quadratic penalty method in~\cite{lin2022complexity} attain $\widetilde{O}(\epsilon^{-3})$ oracle complexity\footnote{We use 
    ${O}(\cdot)$, ${o}(\cdot)$, $ \Omega(\cdot)$, and $\Theta(\cdot)$ to describe asymptotic behaviors. The tilde-notations $\widetilde O(\cdot)$ and $\widetilde \Omega(\cdot)$  additionally suppress polylogarithmic factors. } under suitable regularity conditions (RCs).
    Their total arithmetic complexity is $\widetilde{O}((\mathcal{C}_{\mathrm{FO}} + \mathcal{C}_r)\epsilon^{-3})$.
The inexact quadratically regularized constrained (IQRC) method~\cite{ma2020quadratically} attains $O(\epsilon^{-3})$ oracle and $O((\cC_{\mathrm{FO}} + \cC_r)\epsilon^{-3})$ arithmetic complexity for inequality-constrained smooth problems without a nonsmooth regularizer $r$. 
Allowing a convex regularizer~$r$, the inexact constrained proximal-point method with constraint extrapolation~(ICPPC)~\cite{boob2023stochastic} finds a randomized 
near-$\epsilon$-KKT point with the same oracle and arithmetic complexities as IQRC. Several recent FOMs achieve $O(\epsilon^{-2})$
oracle complexity for nonconvex constrained problems \cite{boob2025level,curtis2025proximal,zhu2026nonconvex} but potentially require higher-order arithmetic complexity. 
The level-constrained proximal gradient (LCPG) method \cite[Corollary~4]{boob2025level} returns a randomized near-$\epsilon$-KKT point with arithmetic complexity of $O(\cC_{\mathrm{FO}}\epsilon^{-2} + (dn + \cC_r)\epsilon^{-4})$.
      For $\mathcal{X}=\RR_+^d$, after computing a normal step~$v^k$, the proximal-gradient method in~\cite{curtis2025proximal} obtains a tangential step~$u^k$ by solving the following linearly constrained nonsmooth subproblem:%
\begin{equation}
    \begin{aligned}
        \min_{u\in\RR^d}\quad
        & \nabla f(x^k)^\top u
        +\frac{1}{2\alpha_k}\|u\|^2
        +\frac{1}{\alpha_k}(v^k)^\top u
        +r(x^k+v^k+u)\\
        \text{s.t.}\quad
        & J_g(x^k)u=0,\qquad x^k+v^k+u\ge 0.
    \end{aligned} \label{subpg}
    \end{equation}
    The main per-iteration arithmetic cost in~\cite{curtis2025proximal} comes from solving~\cref{subpg}. 
    Zhu and Li~\cite{zhu2026nonconvex} propose a smoothed
prox-linear ALM (SPALM) for the following composite constrained problem with smooth $\{c_i\}$  and convex $\{h_i\}$: 
\[\min_{x\in \mathcal X}h_0(c_0(x)),
        \ \mathrm{s.t.}\  
        h_i(c_i(x))\le 0,\quad i=1,\ldots,n.\]
At each outer iteration, SPALM solves a nonsmooth composite subproblem:
    \begin{equation}
    \begin{aligned}
         \min_{x\in\mathcal X}\quad
        & h_0\!\left(c_0(x^k)+\nabla c_0(x^k)^\top(x-x^k)\right)
        +\frac{1}{2\lambda}\|x-x^k\|^2
        +\frac{r_x}{2}\|x-z^k\|^2\\
        &\quad +\frac{\rho}{2}\sum_{i=1}^{n}
        \left(
        \left[\frac{y_i^k}{\rho}+h_i\!\left(c_i(x^k)+\nabla c_i(x^k)^\top(x-x^k)\right)\right]_+^2
        -\left(\frac{y_i^k}{\rho}\right)^2
        \right).
    \end{aligned} \label{sublr}
    \end{equation} 
    SPALM obtains $O(\epsilon^{-3})$ outer-iteration complexity under local conic regularity and $O(\epsilon^{-2})$ under an additional structural condition \cite[Assumption 4.1]{zhu2026nonconvex}. However, neither~\cite{curtis2025proximal} nor~\cite{zhu2026nonconvex} fully accounts for the arithmetic cost of solving~\cref{subpg} or~\cref{sublr}.

To narrow the gap between $O(\epsilon^{-2})$ oracle complexity and 
arithmetic complexity for 
\cref{eq:model}, we introduce a \textbf{S}moothed \mbox{\textbf{E}xact-\textbf{P}enalty} method with \textbf{i}nexact \mbox{\textbf{Prox}-linear} updates (\mbox{SEPi-Prox}). Employing two efficient subsolvers, it retains the $O(\epsilon^{-2})$ oracle complexity and achieves competitive arithmetic complexities; see \cref{tab:complexity_comparison}.

\begin{table}[!htbp]
\centering
\renewcommand{\arraystretch}{1.2}
\resizebox{\textwidth}{!}{
\begin{tabular}{l|c|c|c|c|c}
\toprule
Method
& Eq.$^*$
& Ineq.$^*$
& Regularity Condition 
& \thead{\normalsize Oracle \\ \normalsize Complexity}
& Arithmetic Complexity
\\
\hline

iALM~\cite{li2021rate} 
& \checkmark
& $\times$
& Iterate-wise NSC$^\ddagger$
& $\widetilde{O}(\epsilon^{-3})$
& $\widetilde{O}\!\left(
    (\cC_{\mathrm{FO}}+\cC_r)\epsilon^{-3}
  \right)$
\\
 
iPPP \cite{lin2022complexity} & \checkmark
& \checkmark
& Global NSC$^\ddagger$
& $\widetilde{O}(\epsilon^{-3})$
& $\widetilde{O}\!\left(
    (\cC_{\mathrm{FO}}+\cC_r)\epsilon^{-3}
  \right)$
\\

LCPG \cite{boob2025level} 
& $\times$
& \checkmark
& Uniform MFCQ
& $O(\epsilon^{-2})$
& $O\!\left(
    \cC_{\mathrm{FO}}\epsilon^{-2}
    +(dn+\cC_r)\epsilon^{-4}
  \right)$$^\P$
\\

IQRC$^\dagger$~\cite{ma2020quadratically}/ICPPC~\cite{boob2023stochastic}
& $\times$
& \checkmark
& Uniform Slater/Strong feasibility
& $O(\epsilon^{-3})$
& $O\!\left(
    (\cC_{\mathrm{FO}}+\cC_r)\epsilon^{-3}
  \right)$
\\

Proximal-gradient$^\dagger$ \cite{curtis2025proximal} 
& \checkmark
& $\times$
& Sequential LICQ-type CQ$^\ddagger$
& $O(\epsilon^{-2})$
& $O\!\left(
    (\cC_{\mathrm{FO}}+\cC_{s1})\epsilon^{-2}
  \right)$$^\P$
\\

\hline
\multirow{3}{*}{SPALM \cite{zhu2026nonconvex}} 
& \multirow{3}{*}{$\times$}
& \multirow{3}{*}{\checkmark}
& Local conic regularity$^\ddagger$
& $O(\epsilon^{-3})$
& $O\!\left(
    (\cC_{\mathrm{FO}}+\cC_{s2})\epsilon^{-3}
  \right)$$^\P$
\\
\cline{4-6}
&
& 
&\multirow{2}{*}{   \thead{\normalsize    Local conic regularity$^\ddagger$ \\ \normalsize   + \cite[Assumption 4.1]{zhu2026nonconvex}}}
& 
\multirow{2}{*}{$O(\epsilon^{-2})$}
&
\multirow{2}{*}{$O (
    (\cC_{\mathrm{FO}}+\cC_{s2})\epsilon^{-2}
   )$$^\P$}
\\
&&&&&\\

\hline

Algorithms \ref{alg:spplm}+\ref{alg:acceleratedNesterov} 
& \multirow{2.5}{*}{\checkmark}
& \multirow{2.5}{*}{$\times$}
& \multirow{2.5}{*}{\thead{\normalsize Local NSC \\ \normalsize(\cref{ass:regularity})}}
& \multirow{2.5}{*}{$O(\epsilon^{-2})$}
& $\widetilde{O}\!\left(
    \cC_{\mathrm{FO}}\epsilon^{-2}
    +(dn+\cC_r)\epsilon^{-5/2}
  \right)$
\\

Algorithms \ref{alg:spplm}+\ref{alg:pdmcp}+\ref{alg:vaidya-update}
&  
& 
& 
& 
& $\widetilde{O}\!\left(
    (\cC_{\mathrm{FO}}
    +dn^2+n\cC_r+n^{3.373})\epsilon^{-2}
  \right)$
\\ 
 
\bottomrule
\end{tabular}
}

\vspace{1mm} 
\caption{Comparisons of FOMs for nonlinear nonconvex constrained optimization equipped with a nonsmooth regularizer and a closed convex set.
\textsuperscript{$*$}``Eq.'' and ``Ineq.'' denote  
equality and inequality constraints, respectively.  
\textsuperscript{$\dagger$}For IQRC~\cite{ma2020quadratically}, the complexity bound is established with $r \equiv 0$. The proximal-gradient method~\cite{curtis2025proximal} requires   
$\mathcal X = \RR^d_+$.
\textsuperscript{$\ddagger$}Iterate-wise NSC and Global NSC require the non-singularity condition (NSC, \cref{as.cq}) to hold at all iterates $x^k$ and over the entire set~$\mathcal{X}$, respectively. Sequential LICQ-type CQ and local conic regularity refer to \cite[Assumption~5.1]{curtis2025proximal} and \mbox{\cite[Assumption~2.2]{zhu2026nonconvex}}, respectively.
\textsuperscript{$\P$}The quantities $\cC_{s1}$ and $\cC_{s2}$ denote the arithmetic costs of solving \cref{subpg} and \cref{sublr}, respectively.  
For LCPG, the arithmetic cost corresponds to \cite[Corollary 4]{boob2025level}.} 
\label{tab:complexity_comparison}
\end{table}

\subsection{Technical approach}\label{sec:motivation} 
We apply 
the $\ell_2$-norm function $h(\cdot)=\|\cdot\|$ to the equality constraints and obtain  
the exact-penalty reformulation of \eqref{eq:model} as 
\begin{align}\label{eq:penaltymodel}
	\min_{x\in \mathcal X}\,\, F(x;\beta):=  f(x) + r(x) + \beta h(g(x)), \tag{EP}
\end{align}
where $\beta > 0$ is the penalty parameter. Under appropriate RCs, 
a sufficiently large  
$\beta$ can ensure that any stationary point of the penalized problem~\cref{eq:penaltymodel} yields a KKT point of the original constrained problem~\cref{eq:model}~\cite{byrd2012line}.

Since $h$ is nonsmooth, we adopt a smoothing approach~\cite{nesterov2005smooth,burke2013epi,chen2012smoothing,beck2012smoothing} to facilitate efficient first-order updates. For any smoothing parameter $\nu > 0$, define the smoothed exact-penalty function $h_{\nu}$ as
\begin{align}\label{dualhv}  
h_{\nu}(z) := \max_{y\in \mathcal{B}_1}\left\{\langle y, z\rangle - \frac{\nu}{2}\|y\|^2\right\} = 
\begin{cases}
\frac{\|z\|^2}{2\nu}, & \text{if } \|z\| \le \nu, \\
\|z\| - \frac{\nu}{2}, & \text{otherwise,}
\end{cases}
\end{align}
where $\mathcal{B}_1 := \{y\in \RR^n : \|y\| \le 1\}$.
Smoothed exact-penalty formulations are widely used in constrained optimization, including Huber-type penalties~\cite{nedich2023huber,tatarenko2021smooth,chen2016penalty,liu2026spider} and softplus penalties~\cite{li2023softplus,li2025new}.
We show in Appendix~\ref{apdaux} that
\begin{align}\label{eq:smoothing-variation} 
0 \le h_{\nu_2}(z) - h_{\nu_1}(z) \le \frac{1}{2}(\nu_1 - \nu_2), \quad \forall z\in\mathbb{R}^n, \ \forall 0 < \nu_2 \le \nu_1.
\end{align}
Since $\lim_{\nu_2\to 0} h_{\nu_2}(z) = \|z\|$, it follows from~\eqref{eq:smoothing-variation} that
\begin{align}\label{eq:smt-abs} 
0 \le \|z\| - h_{\nu_1}(z) \le \tfrac{\nu_1}{2}.
\end{align}
In addition, $h_{\nu}$ is $\frac{1}{\nu}$-smooth with $\nabla h_{\nu}(z) = \operatorname{Proj}_{\mathcal{B}_1}(z/\nu)$.
Replacing $h$ with $h_{\nu}$ in~\cref{eq:penaltymodel} yields a smoothed exact-penalty reformulation of \cref{eq:model}:
\begin{align}\tag{SP} \label{eq:sp}
\min_{x\in\mathcal{X}} F_{\nu}(x; \beta) := f(x) + r(x) + \beta h_{\nu}(g(x)).
\end{align}
Inspired by the prox-linear framework in~\cite{drusvyatskiy2019efficiency,zhang2022stochastic}, we perform the following prox-linear updates to \eqref{eq:sp} for each outer iteration $k\ge0$:%
\begin{equation}\label{eq:subp}
\begin{aligned}
\min_{x\in\mathcal{X}} f(x^k) &+ \nabla f(x^k)^{\top}(x - x^k)  + r(x)  \\[-1mm]
&+ \beta h_{\nu}\big(J_g(x^k)(x - x^k)+ g(x^k)\big)+\tfrac{\rho}{2}\|x-x^k\|^2.
\end{aligned}
\end{equation} 
In \Cref{sec:exact}, we analyze how~\cref{eq:sp} and~\cref{eq:subp} relate to  
\cref{eq:model} under a local non-singularity condition (local NSC, \cref{ass:regularity}). 
Under this RC, any near-stationary point $x$ of \cref{eq:sp} is  $\nu$-feasible when $\beta$ is sufficiently large, and the same holds for \cref{eq:subp} when $\|x-x^k\|$ is also sufficiently small. 
At each iteration $k$, our method conducts a backtracking line search for 
$\rho$ and solves the subproblem in \cref{eq:subp} inexactly.  
Solving each subproblem via Nesterov's accelerated proximal-gradient (APG) subsolver~\cite{nesterov2013gradient} or a dual cutting-plane subsolver yields our results in~\cref{tab:complexity_comparison}.

\subsection{Related work}\label{sec:rela}
 
Beyond a detailed comparison to several works in~\cref{tab:complexity_comparison}, we briefly review prior work on FOMs for nonconvex constrained optimization. 

\paragraph{\textbf{Penalty-based methods}}
ALMs~\cite{sahin2019inexact,li2021rate,shi2022momentum,zuo2025adaptive} and quadratic penalty methods~\cite{bourkhissi2025complexity,kong2019complexity,lin2022complexity} are widely used approaches for constrained optimization.
For nonsmooth exact-penalty methods,  existing works~\cite{diouane2026nonsmooth,wang2017penalty,cui2025exact} focus on smooth problems without a nonsmooth regularizer or a convex 
set $\cX$. Smoothed exact-penalty methods~\cite{liu2026spider,xu2026smoothing} have been studied for inequality-constrained nonsmooth problems.

\paragraph{\textbf{SQP and proximal-point methods}} SQP methods solve nonlinearly constrained problems via sequential quadratic subproblems with linearized constraints~\cite{gill2005snopt,berahas2021sequential}, and this paradigm has recently been extended to problems with a nonsmooth regularizer~\cite{dai2025proximal,curtis2025proximal}.
Proximal-point methods address nonconvexity in both the objective and constraints by adding quadratic proximal terms to form a sequence of strongly convex constrained subproblems. Representative approaches include quadratically regularized methods~\cite{ma2020quadratically,boob2023stochastic} and the LCPG method~\cite{boob2025level}.

\paragraph{\textbf{Fletcher's penalty and constraint-dissolving methods}}
The aforementioned methods impose assumptions only up to the first-order derivatives of $f$ and $g$. In contrast, several existing approaches require assumptions on the second-order derivatives of $f$ or $g$.
In particular, Fletcher's penalty methods have been developed for \mbox{equality-constrained} problems with $r\equiv0$ and $\cX=\RR^d$~\cite{cui2026fletcher,fletcher1970class,di1986exact,shi2026curvature}. 
For $r\equiv 0$, Xiao et al.~\cite{xiao2024dissolving} propose a constraint-dissolving approach for equality-constrained optimization on $\cX=\RR^d$, which is extended to general closed convex sets in~\cite{xiao2026exact}.

\subsection{Contributions} 
Our main contributions  
are threefold. First,  
under smoothness conditions in \cref{ass:problemsetup} and the local NSC in \cref{ass:regularity}, which is weaker than several widely used RCs and CQs~\cite{zuo2026cq}, we show that any near-stationary point of \cref{eq:sp} or \cref{eq:subp} is near-feasible for \eqref{eq:model}. 
Second, we develop a first-order smoothed exact-penalty method with inexact prox-linear updates, given in \cref{alg:spplm}, for problem \eqref{eq:model}, where 
a backtracking line search is performed on the proximal parameter. 
We establish \(O(\epsilon^{-2})\) oracle complexity for finding an \(\epsilon\)-KKT point of \eqref{eq:model}, matching the first-order oracle lower bound for unconstrained smooth nonconvex optimization~\cite{carmon2020lower}. Third, we pair \cref{alg:spplm} with two subsolvers for the prox-linear subproblems. By applying Nesterov's APG subsolver, our method attains a total arithmetic complexity of 
$\widetilde O\left(\mathcal C_{\rm FO}\epsilon^{-2}+(dn+\cC_r)\epsilon^{-5/2}\right)$. 
With a dual cutting-plane subsolver based on Vaidya's method, our method achieves $\widetilde O\left((\mathcal C_{\rm FO}+dn^2+n\cC_r+n^{3.373})\epsilon^{-2}\right)$ arithmetic complexity, 
which is favorable when $n$ is small. Numerical experiments on multi-class Neyman--Pearson classification instances demonstrate the effectiveness of our method with both subsolvers. 
	
	\subsection{Notation and definitions}
	For a matrix $A$, we use $\Vert{}A\Vert{}$ to denote its spectral norm.  For a vector $v$, we define $[v]_+ := \max\{v, 0\}$ componentwise. 
    For a nonempty closed set $\mathcal{C} \subseteq \mathbb{R}^d$ and 
    a point $x \in \mathbb{R}^d$, 
    $\mathrm{dist}(x, \mathcal{C}) := \inf_{y \in \mathcal{C}} \|x - y\|$. When $\mathcal{C}$ is additionally convex, $\mathrm{Proj}_{\mathcal{C}}(x) := \arg\min_{y \in \mathcal{C}} \|x - y\|$.  
    The normal cone of $\mathcal{X}$ at $x \in \mathcal{X}$ is defined as $\mathcal{N}_{\mathcal{X}}(x) := \left\{ v \in \mathbb{R}^d : \langle v, z - x \rangle \le 0, \; \forall z \in \mathcal{X} \right\}$.  
    For a proper lower semicontinuous convex function $\psi: \mathbb{R}^d \to \mathbb{R} \cup \{+\infty\}$, its subdifferential at $x \in \mathrm{dom}(\psi)$ is denoted by $\partial \psi(x)$. 
Its proximal mapping is defined by \mbox{$\mathrm{prox}_{\psi}(x) := \mathop{\arg\min}_{z} \left\{  \psi(z) + \frac{1}{2} \|z - x\|^2 \right\}$}.

\begin{definition}[Stationarity]
	\label{def:nearkkt}
	Let $\epsilon\ge0$. 
	  We call $x^*\in \mathcal X$  an \textit{$\epsilon$-KKT point} of problem \eqref{eq:model}, if there exist $\xi\in\partial r(x^*)$ and $\lb\in\RR^n$ such that
	\begin{align*}
		\dist\left(\nabla f(x^*)+\xi +  J_{g}  (x^*) ^{\top}\lb ,-\mathcal N_{\mathcal X}(x^*)\right) \leq \epsilon,\text{ and } \|g(x^*)\|\leq \epsilon.
	\end{align*} 
\end{definition}

	\subsection{Organization} 
\Cref{sec:exact} establishes conditions under which approximate stationarity of the smoothed penalty and prox-linear models implies near-feasibility for problem \cref{eq:model}. \Cref{sec:frmwk} presents the smoothed exact-penalty method with inexact prox-linear updates and establishes its outer-iteration complexity.
\Cref{sec:inner} 
derives the total arithmetic complexity
bounds by specifying two subsolvers.
Numerical experiments are presented in  
\Cref{sec:numerics}.  

	\section{Model analysis}	
	\label{sec:exact} 
    In this section, we relate near-stationary points of the smoothed model \cref{eq:sp} and the prox-linear model \cref{eq:subp}  to near-feasible points of \cref{eq:model}. 
We impose the following RC on the constraint system in \cref{eq:model}.  
    \begin{assumption}[Local NSC]\label{ass:regularity}
There exist constants $\sigma, B>0$ such that
\begin{align}\label{as.cq}
 \sigma\|g(x)\|\le \operatorname{dist}
\bigl(J_g(x)^\top g(x),-\mathcal{N}_{\mathcal X}(x)\bigr) , 
\end{align}
whenever $x\in\mathcal X$ and $\|g(x)\|\le B$. 
    \end{assumption}
    
     \cref{ass:regularity} imposes the local non-singularity condition (NSC) only on the $B$-feasible set $\{x\in \cX:\|g(x)\|\le B\}$. 
    When $\mathcal X$ is compact, both the extended rCRCQ~\cite[Assumption 1.1(2)]{xiao2026exact} and the uniform MFCQ~\cite[Assumption 3]{boob2025level} (after slack reformulation) imply \cref{ass:regularity}. The global NSC for mixed constraints~\cite[Assumption 4]{lin2022complexity} also implies \cref{ass:regularity} after slack reformulation, provided that the Jacobian of the inequality function is bounded on $\mathcal X$. However, none of the converse implications holds in general; see~\cite{zuo2026cq}.

\subsection{Relationship between \cref{eq:model} and \cref{eq:sp}}

The following proposition states that a near-stationary point of \cref{eq:sp} is a near-feasible point of \eqref{eq:model} when the penalty parameter is sufficiently large.

\begin{proposition}\label{thm:rl-model-smooth}
Suppose that \cref{ass:problemsetup,ass:regularity} hold.
Let $\nu>0$, $\epsilon\in(0,1)$, and $x\in\mathcal X$.
Assume that
$\beta\geq (\|\nabla f(x)\|+l_r+1)/\sigma$ and $\|g(x)\|\leq B$.
If there exists $\xi\in\partial r(x)$ such that
\begin{align}\label{eqv}
\operatorname{dist}\left(
\nabla f(x)+\xi+
\beta J_g(x)^\top
\operatorname{Proj}_{\mathcal B_1}
\left(\frac{g(x)}{\nu}\right),
-\mathcal N_{\mathcal X}(x)
\right)
\le\epsilon,
\end{align}
then $\|g(x)\|<\nu$.
\end{proposition}

\begin{proof}  
    By \cref{ass:problemsetup},  $\epsilon<1$, \cref{eqv}, and the conic property of $\mathcal N_{\mathcal X}(x)$, we have \begin{align*}
        &\textstyle \beta \dist\left(J_g(x)^{\top}\mathrm{Proj}_{\mathcal B_1}\left(\frac{g(x)}{\nu}\right),-\mathcal N_{\mathcal X}(x)\right) =  \dist\left(\beta J_g(x)^{\top}\mathrm{Proj}_{\mathcal B_1}\left(\frac{g(x)}{\nu}\right),-\mathcal N_{\mathcal X}(x)\right)\\
         \leq\,& \|-\nabla f(x)-\xi\| + \dist\left(\nabla f(x)+\xi+\beta J_g(x)^{\top}\mathrm{Proj}_{\mathcal B_1}\left(\frac{g(x)}{\nu}\right),-\mathcal N_{\mathcal X}(x)\right)\\ 
         \leq\,& \|-\nabla f(x)-\xi\|+\epsilon \leq \|\nabla f(x)\| + \|\xi\|+\epsilon \leq \|\nabla f(x)\|+l_r+\epsilon,
    \end{align*}
    which implies 
    $\dist(J_g(x)^{\top}\mathrm{Proj}_{\mathcal B_1} (\tfrac{g(x)}{\nu} ),-\mathcal N_{\mathcal X}(x))\leq  (\|\nabla f(x)\|+l_r+\epsilon)/\beta<\sigma$. Suppose, to the contrary, that $\|g(x)\|\ge\nu$. Since $\|g(x)\|\le B$, \cref{as.cq} applies and yields
    {\small\begin{align}\label{gesigma}
\textstyle\dist\left(J_g(x)^{\top}\mathrm{Proj}_{\mathcal B_1} \left(\frac{g(x)}{\nu} \right),-\mathcal N_{\mathcal X}(x)\right)= \dist\left(J_g(x)^{\top} \left(\frac{g(x)}{\|g(x)\|} \right),-\mathcal N_{\mathcal X}(x)\right)\geq \sigma.
    \end{align}}
    This contradicts the strict upper bound above, and hence $\|g(x)\| < \nu$.  
\end{proof}

\subsection{Relationship between \cref{eq:model} and \cref{eq:subp}}
Similarly, a near-stationary point of the linearized model \cref{eq:subp} is 
a near-feasible point of \eqref{eq:model} when the penalty parameter is sufficiently large and $\|x-x^k\|$ is sufficiently small.
 
\begin{proposition}\label{prop:exapen}
Suppose that \cref{ass:problemsetup,ass:regularity} hold. 
Let $\nu >0$, \mbox{$\epsilon\in(0,1)$}, and
$x^k,x\in\mathcal X$.
Assume $
\|g(x)\|\le B$,
$\|x-x^k\|\leq  \min\{\frac{\epsilon}{2\rho},\frac{\sigma}{4 L_g},\sqrt{\sigma\nu /(2l_gL_g)}\}$, and $
\beta\ge
\frac{2(\|\nabla f(x^k)\|+l_r+1)}{\sigma}
$. Let $u := \mathrm{Proj}_{\mathcal B_1}\left(\frac{g(x^k)+J_g(x^k)(x -x^k)}{\nu }\right)$.
If there exists \mbox{$\xi\in\partial r(x)$} such that 
\begin{align}
\operatorname{dist}\left(
\nabla f(x^k)+\xi+
\beta J_g(x^k)^\top
u + \rho(x-x^k),
-\mathcal N_{\mathcal X}(x)
\right)
\le \epsilon/ 2, \label{st:subp}
\end{align} 
then $\|g(x)\|<\nu $.
\end{proposition}

\begin{proof}
    Using \cref{st:subp} and $\|x-x^k\|\le \epsilon/(2\rho)$, we have
    \begin{align*} 
        & \beta \dist\left(J_g(x )^{\top}\mathrm{Proj}_{\mathcal B_1}\left(\tfrac{g(x)}{\nu}\right),-\mathcal N_{\mathcal X}(x )\right) = \dist\left(\beta J_g(x )^{\top}\mathrm{Proj}_{\mathcal B_1}\left(\tfrac{g(x)}{\nu}\right),-\mathcal N_{\mathcal X}(x )\right)\\ 
        \leq\, & \left\|-\nabla f(x^k)-\xi+\beta J_g(x )^{\top}\left(\mathrm{Proj}_{\mathcal B_1}\left(\tfrac{g(x)}{\nu}\right)-u\right)+\beta \left(J_g(x )-J_g(x^k)\right)^{\top}u\right\|\\
        &\hspace{2cm}+\rho\|x-x^k\| + \dist\left( \nabla f(x^k) +\xi + \beta J_g(x^k)^{\top}u+\rho(x-x^k),-\mathcal N_{\mathcal X}(x )\right)\\
        \leq\, &  \|\nabla f(x^k)\| + \|\xi\|+\beta \|J_g(x )\|\left\|\mathrm{Proj}_{\mathcal B_1}\left(\tfrac{g(x)}{\nu}\right)-u\right\|+\beta\|J_g(x )-J_g(x^k)\|\|u\| +\epsilon \\
         \leq \,&\|\nabla f(x^k)\|+l_r+\beta l_g\left\|\tfrac{g(x)-g(x^k)-J_g(x^k) (x-x^k)}{\nu}\right\|+\beta \|J_g(x )-J_g(x^k)\|+\epsilon\\
        \le\,& \|\nabla f(x^k)\|+l_r+\beta l_gL_g\|x-x^k\|^2/(2\nu)+\beta L_g\|x -x^k\|+\epsilon,
    \end{align*} 
where the first equality is from the conic property of
$-\mathcal N_{\mathcal X}(x)$, the first two inequalities are
from the $1$-Lipschitz continuity of the distance function and the triangle inequality,
the third inequality is from the nonexpansiveness of projection,
\cref{ass:problemsetup} and $u\in\mathcal B_1$, and the last is
from the smoothness of $g$.
    Hence,  
    {\small\begin{align*}
     \dist\left(J_g(x )^{\top}\mathrm{Proj}_{\mathcal B_1}\left(\tfrac{g(x)}{\nu}\right),-\mathcal N_{\mathcal X}(x )\right) &\leq \frac{\|\nabla f(x^k)\|+l_r+\epsilon}{\beta} +  L_g\|x -x^k\|+ \frac{l_gL_g\|x-x^k\|^2}{2\nu}\\
    &<\frac{\|\nabla f(x^k)\|+l_r+1}{\beta} + \frac{\sigma}{4}+ \frac{\sigma}{4} 
    \le \sigma.
    \end{align*}}  
    If $\|g(x)\|\ge\nu$, then
  \cref{as.cq} gives the opposite
bound, as in \eqref{gesigma}.
Thus $\|g(x)\|<\nu$. 
\end{proof}

\section{Algorithm framework and outer-iteration complexity}
\label{sec:frmwk}  

In this section, we propose our algorithm SEPi-Prox with backtracking
on the proximal parameter for solving \cref{eq:model}.
At each outer iteration, the algorithm constructs a family of regularized linearized models. It then selects the proximal parameter and next iterate via an inexact stationarity condition paired with a descent acceptance test.

Given $x^k\in\mathcal X$, penalty parameter
$\beta_k>0$, smoothing parameter $\nu_k>0$,  and trial parameter $\rho\ge\underline\rho>0$, consider
a family of subproblems:
\begin{align}\label{eq:sub}
    &\min_{x\in\mathbb{R}^d} \ F^k(x;\rho) :=m^k(x) + \Psi^k(x;\rho), \text{ where}\\
\notag  m^k(x) := \,&f(x^k) + \langle \nabla f(x^k), x-x^k\rangle + \beta_k h_{\nu_k}\big(g(x^k)+J_g(x^k)(x-x^k)\big),\\
\notag   \Psi^k(x;\rho) := \,&r(x) + \delta_{\mathcal{X}}(x) + \frac{\rho}{2}\|x - x^k\|^2.  
\end{align} 
Here, $m^k$ is convex and continuously differentiable, while $\Psi^k(\cdot;\rho)$ and $F^k(\cdot;\rho)$ are $\rho$-strongly convex. 
The parameter $\rho$ controls the proximal regularization
and is determined by backtracking.
By \Cref{ass:problemsetup},
 problem~\eqref{eq:sub} has a unique minimizer. 

For each trial parameter $\rho$, a subsolver computes a trial point $\widetilde x\in\mathcal X$ satisfying
$
\operatorname{dist}
\bigl(0,\partial F^k(\widetilde x;\rho)\bigr)
\le\epsilon_k.
$
The pair $(\widetilde x,\rho)$ is accepted if
$
F_{\nu_k}(\widetilde x;\beta_k)
\le F^k(\widetilde x;\rho).
$
Upon acceptance, we set
$x^{k+1}=\widetilde x$ and $\rho_k=\rho$. Otherwise, we replace $\rho$ by $\gamma \rho$ with $\gamma>1$ and compute a new trial point. 
Lemma~\ref{lm:barrho} established later shows that when $\rho\ge L_f+\beta_kL_g$, the condition $F_{\nu_k}(\widetilde x;\beta_k)\le F^k(\widetilde x;\rho)$ must hold. 
Consequently, each backtracking loop terminates in finitely many trials, provided that the subsolver attains the prescribed stationarity tolerance $\epsilon_k$.

After accepting the new iterate, we update $\beta_{k+1}$ to maintain the penalty lower bound required by \cref{prop:exapen} and initialize the next backtracking
procedure with
$\rho_{\mathrm{start},k+1}=\max\{\rho_{k}/\gamma,\underline\rho\}$. 
The   algorithm framework is given in \cref{alg:spplm}. Two subsolvers will be introduced in Section~\ref{sec:inner} to obtain $\widetilde x$.

\begin{algorithm}[h!]
\caption{\small SEPi-Prox for problem~\cref{eq:model}}\label{alg:spplm}
\begin{algorithmic}[1]

\STATE \textbf{Input:} $x^0\in\mathcal X $ that satisfies $\|g(x^0)\|\le B/2$,  $\rho_{\mathrm{start},0}=\underline{\rho}>0$,$\gamma>1$, and non-increasing positive sequences $\{\nu_k\}_{k\ge0} $ and $\{\epsilon_k\}_{k\ge0} $.\\
\STATE \textbf{Choose} $\beta_0 \ge \max\left\{2(\|\nabla f(x^0)\|+l_r+1)/\sigma,4(f(x^0)+r(x^0)+1/\underline{\rho}-f^*)/B\right\}$.
\FOR{$k = 0,1,\ldots$}  
    \STATE Set $\rho\leftarrow\rho_{\mathrm{start},k}$.
    \WHILE{true} 
    \STATE Apply a subsolver to find an approximate point $\widetilde x\in \mathcal X$ that satisfies  \begin{align}\label{eq:subcond}
 \dist(0,\partial F^k(\widetilde x;\rho))\leq \epsilon_k. 
\end{align} 
    \textbf{if} $F_{\nu_k}(\widetilde x;\beta_k)\le F^k(\widetilde x;\rho)$, \textbf{then break}.  
    \STATE Update $\rho \gets \rho\gamma$.
    \ENDWHILE 
    \STATE Set $x^{k+1} = \widetilde x$, $\rho_{k} = \rho$,  $\rho_{\mathrm{start},k+1} =\max \{\rho_{k}/\gamma,\underline{\rho}\}$ and \[\beta_{k+1}=\max\left\{\beta_k,2(\|\nabla f(x^{k+1})\|+l_r+1)/\sigma\right\}.\]
\ENDFOR
\end{algorithmic}
\end{algorithm}

We first establish a few necessary lemmas. The first one shows that all iterates remain in the $B$-feasible set so that \eqref{as.cq} applies throughout. 
\begin{lemma}\label{lem:bndcon}
Suppose that \Cref{ass:problemsetup} holds.
For all $k\ge0$, we have
\begin{align}\label{idp}
F_{\nu_k}(x^{k+1};\beta_k)
\le
F_{\nu_k}(x^k;\beta_k)
+\frac{\epsilon_k^2}{2\underline\rho}.
\end{align}
Additionally, if $\nu_0\in (0,B/2]$ and  $\sum_{k=0}^\infty\epsilon^2_k \le 2$, then, for all $k\ge0$, it holds that $
h_{\nu_k}(g(x^{k+1}))
\le B-\frac{\nu_k}{2},
$ and $
\|g(x^{k+1})\|\le B.
$
\end{lemma}

\begin{proof}
Choose $\xi^{k+1}\in\partial F^k(x^{k+1};\rho_k)$ with
$\|\xi^{k+1}\|\le\epsilon_k$.
Strong convexity gives
\begin{align*}
F^k(x^{k+1};\rho_k)
&\le
F^k(x^k;\rho_k)
+\langle \xi^{k+1},x^{k+1}-x^k\rangle
-\frac{\rho_k}{2}\|x^{k+1}-x^k\|^2\\
&\le
F^k(x^k;\rho_k)
+\frac{\|\xi^{k+1}\|^2}{2\rho_k}\le
F_{\nu_k}(x^k;\beta_k)
+\frac{\epsilon_k^2}{2\underline\rho},
\end{align*}
where the second inequality is by Young's inequality and the last inequality is because $F^k(x^k;\rho_k)=F_{\nu_k}(x^{k};\beta_k)$ and $\|\xi^{k+1}\|\le\epsilon_k$.
Combining this estimate with the acceptance condition
$
F_{\nu_k}(x^{k+1};\beta_k)
\le F^k(x^{k+1};\rho_k)
$
proves \eqref{idp}. 
Next, we prove $h_{\nu_k}(g(x^{k+1}))\leq B-{\nu_k}/2$ for all $k\geq 0$ by induction.  
Fix $k \ge 0$, and when $k \ge 1$, assume $h_{\nu_i}(g(x^{i+1})) \le B - \nu_i/2$ for $i = 0, \ldots, k - 1$. Summing \eqref{idp} over $i=0,\ldots,k$ gives
\begin{align*}
f(x^{k+1})+r(x^{k+1})
+\sum_{i=0}^{k}\beta_i h_{\nu_i}(g(x^{i+1}))&\le
f(x^0)+r(x^0)
+\sum_{i=0}^{k}\beta_i h_{\nu_i}(g(x^i))
+\sum_{i=0}^{k}\frac{\epsilon_i^2}{2\underline\rho}\\
&\le
f(x^0)+r(x^0)
+\sum_{i=0}^{k}\beta_i h_{\nu_i}(g(x^i))
+\frac{1}{\underline\rho}.
\end{align*} 
Rearranging terms of the above inequality gives 
\begin{equation}\label{eq:sum-beta-h}
\begin{aligned} 
    &\beta_k h_{\nu_k}(g(x^{k+1})) +f(x^{k+1}) + r(x^{k+1}) -f(x^{0}) - r(x^{0})-1/\underline{\rho} \\  
      \le  \,& \beta_0 h_{\nu_0}(g(x^{0}))+  \sum_{i=0}^{k-1} \big(\beta_{i+1}h_{\nu_{i+1}}(g(x^{i+1}))-\beta_i h_{\nu_i}(g(x^{i+1}))\big). 
\end{aligned}  
\end{equation}
For the second term in RHS of \eqref{eq:sum-beta-h}, we have  
\begin{align}
    &   \sum_{i=0}^{k-1} \left(\beta_{i+1}h_{\nu_{i+1}}(g(x^{i+1}))-\beta_i h_{\nu_i}(g(x^{i+1}))\right) \label{fupb} \\
     \le   &\sum_{i=0}^{k-1} \left(\frac{\beta_{i+1}(\nu_i-\nu_{i+1})}2+(\beta_{i+1}-\beta_i)\left( B -  \frac{\nu_i} 2  \right)\right)   
     =   \frac{\beta_{0}\nu_0 -\beta_{k}\nu_{k}}2+(\beta_{k}-\beta_0)B, \notag
\end{align} 
where both sides vanish when $k = 0$. For $k \ge 1$, the inequality follows from \cref{eq:smoothing-variation}, the induction hypothesis, and
$\beta_{i+1}\ge\beta_i$. 
By \cref{eq:smt-abs} and the initialization, \mbox{$h_{\nu_0}(g(x^{0}))\le \|g(x^{0})\|\le B/2$} and $\nu_0 \le B/2$. Combining \eqref{eq:sum-beta-h} and \eqref{fupb} then gives
\begin{equation}\label{fupb2}
\begin{aligned}
    &\beta_k h_{\nu_k}(g(x^{k+1})) +f(x^{k+1}) + r(x^{k+1}) -f(x^{0}) - r(x^{0})-1/\underline{\rho} \\
     \le\,& \beta_0 h_{\nu_0}(g(x^{0}))+\beta_0\nu_0/2-\beta_0 B+  \beta_{k}\left( B -  \nu_k/2   \right)\\
     \le\,&  \beta_0 B/2+  \beta_0\nu_0/2-\beta_0 B+  \beta_{k}\left( B -  \nu_k/2   \right)  \le   -\beta_0 B/4+\beta_{k} \left( B -  \nu_k/2   \right)  .
\end{aligned}
\end{equation}
Since $f(x^{k+1})\ge f^*$, $r(x^{k+1})\ge 0$ and $\beta_0 \geq {4(f(x^0)+r(x^0)+1/\underline{\rho}-f^*) }/{B}$, \cref{fupb2} implies  $h_{\nu_k}(g(x^{k+1}))\le B-\nu_k/2,$
which completes the induction.  
Finally, from \cref{eq:smt-abs}, we have 
\(\|g(x^{k+1})\| \leq h_{\nu_k}(g(x^{k+1}))+\nu_k/2\le B\) for all $k\ge0$.
 \end{proof}

The lemma below shows that the nondecreasing sequence $\{\beta_k\}_{k\ge 0}$ is bounded.

\begin{lemma}\label{lem:barbeta}
Suppose that \Cref{ass:problemsetup} holds,
$
\sum_{i=0}^{\infty}\epsilon_i^2\le2,
$ and $
\nu_0\in(0,B/2].$
Then 
$
\beta_k\le\bar\beta:=\max\{\beta_0,\widehat\beta\},
$ and $
\|\nabla f(x^k)\|\le l_f:=
\max\left\{
\|\nabla f(x^0)\|,
\sqrt{2L_f\bar\beta B}
\right\}
$ for all $k\ge0$,
where  
\begin{align}\label{eq:hat-beta}
\widehat\beta
:=
\frac{
2(l_r+1)\sigma+4L_fB
+4\sqrt{L_fB(l_r+1)\sigma+L_f^2B^2}
}{\sigma^2}.
\end{align}
\end{lemma}

 \begin{proof}
It is trivial that $\beta_0\le\bar \beta$.
Assume $\beta_k\leq \bar \beta$ for some $k\ge 0$.
Then 
\begin{align*}
    f(x^{k+1}) - f^*  
    &\le f(x^{0}) + r(x^{0}) -f^*+1/\underline{\rho}-\beta_0 B/4+ \beta_{k}\left( B - \nu_k/2   \right)  \\
    &\le \beta_k(B-\nu_k/2) \leq \beta_kB\leq \bar \beta B,
\end{align*}
where the first inequality is from $h_{\nu_k}(g(x^{k+1}))\ge 0$, $r(x^{k+1})\ge0$ and \cref{fupb2},  
and the second one is by the choice of $\beta_0$.  

From the $L_f$-smoothness of $f$ on $\RR^d$, we have
$
    \|\nabla f(x^{k+1})\|^2 \leq 2L_f (f(x^{k+1})-f^*).
$
 The bound on $f(x^{k+1})-f^*$ obtained above therefore gives
\mbox{$
\|\nabla f(x^{k+1})\|
\le\sqrt{2L_f\bar\beta B}.
$}
Hence, once $\beta_k \le \bar{\beta}$ holds for all $k \ge 0$, it follows that \mbox{$\|\nabla f(x^k)\| \le l_f$} for all $k \ge 0$.

To verify that the   update preserves the bound
$\bar\beta$, let $u>0$ and   $t=\sqrt u$. Then
\begin{align*}
& {2\left(\sqrt{2L_fBu}+l_r+1\right)}/{\sigma}\le u
\iff
\sigma t^2-2\sqrt{2L_fB}\,t-2(l_r+1)\ge 0 \iff
u\ge
\widehat{\beta}.
\end{align*}
Since $\bar{\beta}\ge \widehat{\beta}$, we have
\(
 {2(\sqrt{2L_fB\bar{\beta}}+l_r+1)}/{\sigma}
\le \bar{\beta}.
\)
Then, 
by \mbox{$\|\nabla f(x^{k+1})\|\leq \sqrt{2L_f\bar \beta B}$}, \mbox{$\beta_k\le \bar \beta$}, and the updating rule of $\beta_{k+1}$, we have  
\begin{align*}
    \beta_{k+1} 
    \le \max\left\{2\left( \sqrt{2L_f\bar \beta B}+l_r+1\right)/\sigma,\beta_k\right\}\le\max\{\bar \beta,\beta_k\} = \bar \beta. 
\end{align*}    
 \end{proof}

By Lemma~\ref{lem:barbeta}, if $\beta_0$ is large enough, then $\{\beta_k\}_{k\ge 0}$ will be a constant sequence.

\begin{corollary}\label{sfclgb}
Under the conditions of \cref{lem:barbeta}, suppose that $\beta_0\ge\widehat\beta$, where $\widehat\beta$ is defined
in \eqref{eq:hat-beta}.
Then $\{\beta_k\}_{k\ge 0}$   generated by
\cref{alg:spplm} satisfies
\mbox{$\beta_k =\beta_0, \forall\, k\ge0.$}
\end{corollary}

\begin{proof}
The additional assumption gives $\bar\beta=\beta_0$.
By \cref{lem:barbeta}, $\beta_k\le\beta_0$.
The update rule also gives $\beta_k\ge\beta_0$.
Hence $\beta_k=\beta_0$ for all $k\ge0$.
\end{proof}

The next lemma shows that the proximal parameter $\rho$ remains bounded during the backtracking line search.

\begin{lemma}\label{lm:barrho}
Under the conditions of \cref{lem:barbeta}, every trial parameter $\rho$ generated by the line search in
\cref{alg:spplm} satisfies
$
\underline\rho\le\rho\le\bar\rho:=
\max\{\gamma(L_f+\bar\beta L_g),\underline\rho\}.
$
In particular, $\rho_k\le\bar\rho$ for all $k\ge0$.
\end{lemma}

\begin{proof}
    From the smoothness of $f$ and $g$,  
it holds for all $x\in \mathcal X$ that   
    \begin{align*}
        &f(x)\le f(x^k)+\langle \nabla f(x^k),x-x^k\rangle +\frac{L_f}{2}\|x-x^k\|^2,\\
        &\|g(x)-g(x^k)-J_g(x^k)(x-x^k)\| \le \frac{L_g}{2}\|x-x^k\|^2. 
    \end{align*}
    Since $h_{\nu_k}$ is 1-Lipschitz continuous,  
    we obtain 
    \begin{align*}
        &F_{\nu_k}(x;\beta_k)-r(x)= f(x)+\beta_k h_{\nu_k}(g(x))\\
        \le & f(x^k)+\langle \nabla f(x^k),x-x^k\rangle + \beta_k h_{\nu_k}(g(x^k)+J_g(x^k)(x-x^k)) + \tfrac{L_f+\beta_kL_g}{2}\|x-x^k\|^2.
    \end{align*}
 Therefore, every trial with
$\rho\ge L_f+\beta_kL_g$ is accepted. 
For any $k$, if $\rho_k>\rho_{\mathrm{start},k}$, $\rho=\rho_k/\gamma$ must be a rejected trial, which means 
$\rho_k/\gamma<L_f+\beta_kL_g$
and thus \mbox{$\rho_k\leq \gamma(L_f+\bar\beta L_g)
\le\bar\rho$}. Therefore, $\rho_k\leq\max\{\rho_{\mathrm{start},k},\bar\rho\}$ for all $k$. 
By this result and Line~9 of \cref{alg:spplm},  it holds that
\[\underline\rho\leq\rho_{\mathrm{start},k+1}=
\max\{\rho_k/\gamma,\underline\rho\}
\le\max\{\rho_k,\underline\rho\}\le\max\{\rho_{\mathrm{start},k},\bar\rho\}.\] Since $\rho_{\mathrm{start},0}=\underline\rho$,
applying this recursively gives $\rho_{\mathrm{start},k}\leq \bar\rho$ for all $k$. This also implies $\rho_k\leq\max\{\rho_{\mathrm{start},k},\bar\rho\}\leq\bar\rho$. Moreover, $\rho_k\geq\rho_{\mathrm{start},k}\geq\underline\rho$ and every trial parameter
$\rho$ satisfies $\underline\rho\le\rho\le\bar\rho$ because it is between $\rho_{\mathrm{start},k}$ and $\rho_k$.
\end{proof}

Now we are ready to  
establish the outer-iteration complexity bound of Algorithm~\ref{alg:spplm} for obtaining an $\epsilon$-KKT point of problem \eqref{eq:model}.

\begin{theorem}
\label{thm:subconv}
Suppose that \cref{ass:problemsetup,ass:regularity} hold.
Fix $\epsilon\in(0,1)$, and choose
 $ \nu_k = \frac{c_\nu}{\sqrt{k+2}\ln(k+2)} $ with $c_\nu\in(0,\frac{B \sqrt{2}\ln 2}2]$, and $
\sum_{i=0}^{\infty}\epsilon_i^2\le2.
$ Then \cref{alg:spplm} must find an $\epsilon$-KKT point of problem \cref{eq:model} 
      within $K$ outer iterations, where   
\begin{align}\label{eqK}
    K =\left\lceil \max\left\{\frac{4M_1(1+\bar \rho +L_f+ \bar \beta  L_g)^2}{ \epsilon^{2}},\frac{16M_1 L_g^2}{ \sigma^2 },6+32M_2^2\ln^2\!\bigl(4M_2+1\bigr)\right\}\right\rceil,
\end{align} 
    with $M_1 := \frac{4 \bar \beta B }{\underline{\rho} }  +  \frac{2(4+\underline{\rho}^2)}{\underline{\rho}^2 }+3c_\nu^2$ and \mbox{$M_2 := 2M_1 l_g L_g/(\sigma c_\nu).$}
\end{theorem}

\begin{proof}

By \cref{eq:subcond}, there exists 
$v^{k+1}\in \partial F^k(x^{k+1};\rho_k)$ such that 
$\|v^{k+1}\|\le \epsilon_k$. Hence, from the $\rho_k$-strong convexity of $F^k(\cdot;\rho_k)$, it follows that
\begin{equation}\label{eq:Fk-dec-step}
\begin{aligned}
F^k(x^{k+1};\rho_k)
&\le F^k(x^k;\rho_k)
 + \langle v^{k+1},x^{k+1}-x^k\rangle
 - \frac{\rho_k}{2}\|x^{k+1}-x^k\|^2 \\
&\le F_{\nu_k}(x^{k};\beta_k)
 + \epsilon_k\|x^{k+1}-x^k\|
 - \frac{\rho_k}{2}\|x^{k+1}-x^k\|^2  \\
&\le  
 F_{\nu_k}(x^{k};\beta_k)-\frac{\rho_k}{4}\|x^{k+1}-x^k\|^2+\frac{\epsilon^2_k}{\rho_k},
\end{aligned} 
\end{equation}
where the second inequality is by $F^k(x^k;\rho_k)=F_{\nu_k}(x^{k};\beta_k)$ and the Cauchy–Schwarz inequality, and the third one is from Young's inequality.

Together with $F_{\nu_k} (x^{k+1};\beta_k)\leq F^k(x^{k+1};\rho_k)$ and $ \rho_k\ge\underline{\rho}$, \cref{eq:Fk-dec-step} implies 
\begin{equation}\label{Fxkxk1}
    \begin{aligned} 
    \|x^{k+1}-x^k\|^2&\le \frac{4}{\rho_k}\left(F_{\nu_k}(x^{k};\beta_k)-F_{\nu_k} (x^{k+1};\beta_k)+\frac{\epsilon^2_k}{\rho_k}\right )  \\
    &\leq \frac{4}{\underline{\rho} }\left(F_{\nu_k}(x^{k};\beta_k)-F_{\nu_k} (x^{k+1};\beta_k)+\frac{\epsilon^2_k}{\underline{\rho} }\right ).
    \end{aligned}
\end{equation} 
Moreover, by adding the nonnegative term $\beta_{K-1}h_{\nu_{K-1}}(g(x^{K}))$, we obtain 
\begin{align*}
    &\sum_{k=0}^{K-1} \left(F_{\nu_k}(x^k;\beta_k)-F_{\nu_k}(x^{k+1};\beta_k)\right)\\
     \le\,&\sum_{k=0}^{K-1} \left(F_{\nu_k}(x^k;\beta_k)-F_{\nu_k}(x^{k+1};\beta_k)\right) +\beta_{K-1}h_{\nu_{K-1}}(g(x^{K})) \\
     \le\,& f(x^{0}) + r(x^{0}) - f^*+\beta_0 h_{\nu_0}(g(x^{0}))+ \beta_0\nu_0/2-\beta_{K-1}\nu_{K-1}/2+(\beta_{K-1}-\beta_{0})B\\
     \le\,& f(x^{0}) + r(x^{0}) - f^*-\beta_0 B/4+\beta_{K-1}(B-\nu_{K-1}/2)\\
     \le \,&\beta_{K-1}(B-\nu_{K-1}/2) \le \beta_{K-1}B \le \bar \beta B,
\end{align*}
where the second inequality is from $f(x^{K})\ge f^*$, $r(x^{K})\ge0$, and \cref{fupb}, the third one holds by $\nu_0 = c_\nu/(\sqrt{2}\ln 2)\le B/2$ and $h_{\nu_0}(g(x^{0}))\le \|g(x^{0})\|\le B/2$ from \cref{eq:smt-abs}, and the fourth one is due to the choice of $\beta_0$.
Hence, 
by \cref{Fxkxk1} and the fact that $\sum_{i=0}^\infty\epsilon^2_i\le 2$, we obtain 
\begin{align*}
&\min_{k=0,\ldots,K-1} \left(\|x^{k+1}-x^k\|^2+\epsilon_k^2+\nu_k^2\right) \le \frac{1}{K}\sum_{k=0}^{K-1} \left(\|x^{k+1}-x^k\|^2 +\epsilon_k^2+\nu_k^2\right) \\
    & \leq \frac{1}{ K} \sum_{k=0}^{K-1} \left(\frac{4}{\underline{\rho} }(F_{\nu_k}(x^k;\beta_k)-F_{\nu_k}(x^{k+1};\beta_k))+ \frac{(4+\underline{\rho}^2) \epsilon_k^2}{\underline{\rho}^2}  +  \frac{c_\nu^2}{(k+2)\ln^2(k+2)}\right) \\ 
    &\le  \frac{4 \bar \beta B }{\underline{\rho} K}  + \frac{1}{K}\left(\frac{2(4+\underline{\rho}^2)}{\underline{\rho}^2 }+3c_\nu^2\right) 
     = \frac{M_1}{  K}\leq  \min\left\{\frac{\epsilon^2}{4(1+\bar \rho +L_f+ \bar \beta  L_g )^2},\frac{\sigma^2  }{16 L_g^2 }\right\},
\end{align*}
where we have used $\sum_{k= 0}^{\infty} \frac{1}{(k+2)\ln^2(k+2)} \le 3$ in the third inequality, and the last inequality is from the first and second lower bounds of $K$ in \cref{eqK}. Furthermore, by the third lower bound on $K$ in \cref{eqK}, it follows from \cref{lem:dxnu} with $u=K$ and $a= M_2$ that  $ K/(\sqrt{K+2}\ln(K+2)) \ge M_2= 2 M_1 l_g L_g/(\sigma c_\nu).$
Thus it holds that
\[
\frac{M_1}{K}
\le \frac{\sigma c_\nu}{2l_gL_g\sqrt{K+2}\ln(K+2)} =\frac{\sigma\nu_K}{2l_gL_g}
\le \frac{\sigma\nu_k}{2l_gL_g},\quad \forall\, 0\le k\le K.
\]

Choose $k \in \argmin _{t=0,\ldots,K-1} \left(\|x^{t+1}-x^t\|^2+\epsilon_t^2+\nu_t^2\right)$. Since $\rho_k\le \bar{\rho}$, we have 
\begin{align}\label{iterbd}
    \max\{\|x^{k+1}-x^k\|,\epsilon_k,\nu_k\}\le \min\left\{\frac{\epsilon}{2(1+\rho_k+L_f+ \bar \beta  L_g )},\frac{\sigma}{4L_g},\sqrt{\frac{\sigma \nu_k}{2l_gL_g}}\right\}.
\end{align} 
Let ${\lambda^k} = \beta_k\mathrm{Proj}_{\mathcal B_1}(\frac{g(x^{k})+J_g(x^k)(x^{k+1}-x^{k})}{\nu_k}).$ By \eqref{eq:subcond} and \cref{iterbd}, there exists a subgradient \mbox{$\xi^{k+1} \in \partial r(x^{k+1})$} such that
\begin{equation}\label{eq:bd-stat-fk}
    \begin{aligned}
    &\dist(\nabla f(x^k) +\xi^{k+1} +  J_g(x^k)^{\top}\lambda^k + \rho_k(x^{k+1}-x^k),-\mathcal N_{\mathcal X}(x^{k+1})) \\
     =\,& \dist(0,\partial F^k(x^{k+1};\rho_k)) \le \epsilon_k\le \epsilon/2. 
    \end{aligned}
\end{equation}  
By Lemma~\ref{lem:bndcon}, we have $\|g(x^{k+1})\|\le B$. Hence, 
by \cref{iterbd} and \cref{prop:exapen}  
with
$x=x^{k+1}$, $\beta = \beta_k$, and $\nu=\nu_k$, we have  
\(
\|g(x^{k+1})\|<\nu_k\le \epsilon/2.
\)
Notice $\|\lambda^k\|\leq \beta_k \le \bar \beta$. Thus from \cref{iterbd}, \cref{eq:bd-stat-fk}, and the Lipschitz continuity of $\nabla f$ and $J_g$, we obtain 
\begin{align*}
    &\dist(\nabla f(x^{k+1}) +\xi^{k+1} +  J_g(x^{k+1})^{\top}{\lambda^k},-\mathcal N_{\mathcal X}(x^{k+1})) \\
     \le\,&\epsilon_k+(L_f+L_g\|\lambda^k\|+ \rho_k)\|x^{k+1}-x^k\|\le\epsilon/2+(L_f+\bar \beta L_g+  \rho_k)\|x^{k+1}-x^k\|\leq  \epsilon. 
\end{align*}  
\end{proof}

\begin{remark}\label{rm:oraclecmp}
At each outer iteration of \cref{alg:spplm}, since \mbox{$\underline{\rho} \le \rho_{\mathrm{start},k} \le \rho_k \le \bar{\rho}$}, the subsolver is invoked at most $\lceil \log_\gamma(\bar{\rho}/\underline{\rho}) \rceil + 1 = O(1)$ times. Therefore, \cref{alg:spplm} requires at most $( \lceil \log_\gamma(\bar{\rho}/\underline{\rho}) \rceil+1) K=O(\epsilon^{-2})$ subsolver calls and oracle evaluations to obtain an \mbox{$\epsilon$-KKT} point. This bound matches the $\Omega(\epsilon^{-2})$ first-order oracle lower bound for smooth nonconvex optimization~\cite{carmon2020lower}.
\end{remark}

\section{Two subsolvers and total arithmetic complexity}\label{sec:inner}

The \(O(\epsilon^{-2})\) outer-iteration complexity in \cref{thm:subconv} does not include
the cost of solving the subproblem \cref{eq:sub}.  
In this section, we develop two subsolvers and analyze the total arithmetic complexity. We first apply Nesterov's APG method \cite{nesterov2013gradient}  to each subproblem. Then we give a dual cutting-plane method, which reduces the total arithmetic complexity when the number of constraints is relatively small. 

\subsection{Nesterov's APG subsolver} 

At iteration $k$, for any $\rho>0$ in \cref{eq:sub},
 $m^k$  is convex and continuously differentiable, while 
$\Psi^k(\cdot;\rho)$ is $\rho$-strongly convex. Evaluating the proximal mapping of $\Psi^k(\cdot;\rho)$
mapping reduces to a single call of
$\operatorname{prox}_{c r+\delta_{\mathcal X}}$ for a scalar $c>0$. We therefore apply Nesterov's APG method
in \cref{alg:acceleratedNesterov}. The next lemma gives the iteration complexity of Nesterov's APG method to find an approximate point that satisfies \cref{eq:subcond} for each  $\rho\ge \underline{\rho}$. Its proof is provided in Appendix~\ref{apdx:apg}.

\begin{lemma}\label{lem:T}
    Suppose that \cref{ass:problemsetup} and the conditions of
\cref{thm:subconv} hold. Apply \cref{alg:acceleratedNesterov} to \cref{eq:sub} with
$ \widehat f(\cdot) = m^k(\cdot), \widehat\Psi(\cdot) = \Psi^k(\cdot;\rho)$,  $ x^0=x^k$, $\mu=\rho\ge \underline{\rho}$, and tolerance~$\epsilon=\epsilon_k$. Then the method returns a point
$\widetilde x$ satisfying \cref{eq:subcond} within
    {\small \begin{align*}
        T_1^k = \left\lceil  
\max\left\{
    \frac{1}{\ln 2},\;
    2\sqrt{\frac{2\gamma_u \widetilde{L}_k}{ \underline{\rho}}}
\right\}
\left[\ln\left(
    \frac{(1+\gamma_u) \left(l_f+l_r+\bar \beta l_g \right)  \widetilde{L}_k \sqrt{ \frac{2\gamma_u\widetilde L_k }{\underline{\rho}}} }{\underline{\rho} \epsilon_k }
\right)\right]_+\right\rceil
+ 2 
    \end{align*} }
    iterations, where  $\widetilde{L}_k :=  \bar \beta l_g^2/{\nu_k}$ and $\gamma_u>1$ is the input parameter of \cref{alg:acceleratedNesterov}. 
\end{lemma}

Then the following theorem gives the total arithmetic complexity using Nesterov's APG method as the subsolver to find an $\epsilon$-KKT point of problem \cref{eq:model}.

\begin{theorem} \label{thm:cmpapg}
Suppose that \cref{ass:problemsetup,ass:regularity} hold. Let $\epsilon \in (0,1)$, and set
\begin{align}\label{thmcdt}
    \nu_k = \frac{c_\nu}{\sqrt{k+2}\ln(k+2)} \quad
 \text{and}\quad 
\epsilon_k=\frac{\epsilon_0}{k+1},
\end{align} 
where $c_\nu\in(0,B\sqrt{2}\ln 2/2]$ and $\epsilon_0\in(0,1]$.
Then \cref{alg:spplm}, with each prox-linear subproblem solved by
\cref{alg:acceleratedNesterov}, finds an $\epsilon$-KKT point of problem \cref{eq:model} with a total arithmetic complexity of
\(
\widetilde{O}\left( \mathcal C_{\mathrm{FO}}\epsilon^{-2} +(dn+\mathcal C_{r})\,\epsilon^{-5/2} \right).
\)
\end{theorem}
\begin{proof} Since $\sum_{i=0}^{\infty}\epsilon_i^2\le2$, by Theorem~\ref{thm:subconv}, \cref{alg:spplm} finds an $\epsilon$-KKT point of problem \cref{eq:model}  within $K$ outer iterations with $K$ defined in \eqref{eqK}. By \cite[Lemma 6]{nesterov2013gradient}, the total number of
inner line search iterations in \cref{alg:acceleratedNesterov}  over   
at most $T_1^k$ APG iterations is bounded by
\[
N_1^k
:=
\left(1+\frac{\ln\gamma_d}{\ln\gamma_u}\right)T_1^k
+\frac{1}{\ln\gamma_u}
\left[\ln\frac{\gamma_u\widetilde L_k}{\gamma_d L_0}\right]_+,
\]
where $\gamma_d,\gamma_u$ and $L_0$ are user-specified positive constants in \cref{alg:acceleratedNesterov}. 
Since \mbox{$K=O(\epsilon^{-2})$}, for all $k\leq K$, we have 
$ \nu_k^{-1} = O(\epsilon^{-1}\ln(1/\epsilon)) $ and $ \epsilon_k^{-1}=O(\epsilon^{-2}). $
Therefore, by \cref{lem:T}, $\widetilde L_k =O(\nu_k^{-1})= \widetilde O(\epsilon^{-1})$ and \mbox{$T_{1}^k=\widetilde O(\widetilde L_k^{1/2})=\widetilde O(\epsilon^{-1/2})$}, which further implies $N_{1}^k=\widetilde O(\epsilon^{-1/2})$ for all $k\leq K$.
Since the line search in \cref{alg:spplm} invokes the subsolver only $O(1)$ times per outer iteration due to the bounds
$\underline\rho\le\rho\le\bar\rho$, the total number of inner APG line search iterations is \( {O}(KN_1^K) = \widetilde{O}(\epsilon^{-5/2}). \) Each APG iteration requires $O(dn+\mathcal C_{r})$ arithmetic operations. Hence, the total arithmetic complexity is \( O(\mathcal C_{\mathrm{FO}}K+(dn+\mathcal C_r)KN_1^K) = \widetilde{O}( \mathcal C_{\mathrm{FO}}\epsilon^{-2} +(dn+\mathcal C_{r})\,\epsilon^{-5/2}). \)
\end{proof}

\subsection{Dual cutting-plane subsolver}

We next derive a separation oracle for the dual of
\cref{eq:sub}. Throughout this subsection, we consider a fixed $k$ and a fixed
trial parameter $\rho\ge\underline\rho$. 

Set
$
b_k:=g(x^k)-J_g(x^k)x^k
$
and
$
\widetilde f^k(x)
:=
f(x^k)+\langle\nabla f(x^k),x-x^k\rangle
+\frac{\rho}{2}\|x-x^k\|^2.
$
Using  \eqref{dualhv}, we   rewrite subproblem \eqref{eq:sub} as the saddle-point problem
\[
\min_{x\in\mathbb R^d}\max_{y\in\mathcal B_1}\phi_k(x,y)
:=
\widetilde f^k(x)+r(x)+\delta_{\mathcal X}(x)
+\beta_k\langle y,J_g(x^k)x+b_k\rangle
-\frac{\beta_k\nu_k}{2}\|y\|^2.
\]
For each $y\in\mathbb R^n$, the function
$\phi_k(\cdot,y)$ is proper closed and
$\rho$-strongly convex. It therefore has a unique minimizer, denoted as 
\[
x(y)
= \argmin_{x} \phi_k(x,y) = 
\operatorname{prox}_{r/\rho+\delta_{\mathcal X}}
\left(
x^k-\frac{\nabla f(x^k)+\beta_kJ_g(x^k)^\top y}{\rho}
\right).
\]
Define
$
D_k(y):=\phi_k(x(y),y).
$
The dual problem is $\max_{y\in\mathcal B_1}D_k(y)$.

 By Danskin's Theorem \cite{bertsekas1999nonlinear}, $D_k(y)$ is differentiable, and for all $y\in \mathcal B_1$, we have  
 \[\nabla D_k(y) = \beta_k(J_g(x^k) x(y)+b_k-{\nu_k} y).\] 
Moreover,
$D_k$ is
$\beta_k\nu_k$-strongly concave and admits a unique maximizer
$\bar y$ on $\mathcal B_1$.

The optimality condition for maximizing \(D_k\) over
\(\mathcal B_1\) is
$
0\in -\nabla D_k(\bar y)+\mathcal N_{\mathcal B_1}(\bar y).
$
Equivalently,
$
0\in -\beta_k(J_g(x^k)x(\bar y)+b_k-{\nu_k}\bar y)+\mathcal N_{\mathcal B_1 }(\bar y),
$ and hence
\begin{align}\label{eq:pdopt}
\bar y=\operatorname{Proj}_{\mathcal B_1}
\left(\tfrac{J_g(x^k)x(\bar y)+b_k}{{\nu_k}}\right). 
\end{align} 
Define 
$
    R(y):=\beta_k \|J_g(x^k) ^{\top}
(\mathrm{Proj}_{\mathcal B_1}(\tfrac{J_g(x^k)  x(y)+b_k}{{\nu_k}})-y)\|.
$
Since $x(y)$ minimizes~$\phi_k(\cdot,y)$, we have $0\in \partial_x\phi_k(x(y),y)$, and $R(y)$  bounds the subproblem stationarity residual:
{ 
\begin{equation}
\begin{aligned}
   & \operatorname{dist}\left( 0,\partial F^k(x(y);\rho)\right) \\
    =\,&
\operatorname{dist}\left( 0,\partial_{ x}\phi_k(x(y),y)
+
\beta_k J_g(x^k)^{\top}
\left(\mathrm{Proj}_{\mathcal B_1}\left(\tfrac{J_g(x^k)  x(y)+b_k}{{\nu_k}}\right)-y\right) \right) \\
 \le\,&
\operatorname{dist}\left( 0,\partial_{ x}\phi_k(x(y), y)\right)
+
\beta_k \left\|J_g(x^k)^\top\left(\mathrm{Proj}_{\mathcal B_1}\left(\tfrac{J_g(x^k)  x(y)+b_k}{{\nu_k}}\right)- y\right)\right\|= R(y). \label{boundbyR} 
\end{aligned} 
\end{equation}}

We then present an upper bound for $\|\nabla D_k( y)\|$ on $\cB_1$ in the next lemma. 
\begin{lemma} \label{lem:Ck}
Suppose the conditions of \cref{lem:barbeta} hold.
For any trial parameter $\rho\ge\underline\rho$, we have \[\max_{ y\in\mathcal B_1}\|\nabla D_k( y)\|\leq C_k:= \beta_k l_g \frac{l_f + \beta_k l_g+l_r}{\underline{\rho}}+\beta_k  B+\beta_k{\nu_k} ,\]
\end{lemma}
where $l_f$ is defined in \cref{lem:barbeta}.

\begin{proof}
Since $x^k\in\mathcal X$, we have
$0\in \partial \delta_{\mathcal X}(x^k) \ =  \mathcal N_{\mathcal X}(x^k)$. Fix $y\in\mathcal B_1$ and choose
$\xi^k\in\partial r(x^k)$. Then
\(
\zeta_y
:=
\nabla f(x^k)+\beta_kJ_g(x^k)^\top y+\xi^k
\in\partial_x\phi_k(x^k,y),
\)
and
\mbox{$
\|\zeta_y\|\le l_f+\beta_k l_g+l_r.
$}
Since $0\in\partial_x\phi_k(x(y),y)$, strong monotonicity gives
\[
\rho\|x(y)-x^k\|^2
\le
\langle\zeta_y,x^k-x(y)\rangle
\le
(l_f+\beta_k l_g+l_r)\|x(y)-x^k\|.
\]
Therefore,
$
\|x(y)-x^k\|
\le {(l_f+\beta_k l_g+l_r)}/{\rho}.
$
Using \cref{lem:bndcon}, we obtain $\|g(x^k)\|\le B$ and 
\begin{align}
\|\nabla D_k( y)\| &= \beta_k\|g(x^k) + J_g(x^k)(x(y)-x^k)-{\nu_k} y\|\notag\\
        &\le \beta_k\|J_g(x^k) (x(y)-x^k)\| +\beta_k \|g(x^k)\|+\beta_k{\nu_k}\|y\|\notag\\
&\le
\beta_k l_g\|x(y)-x^k\|
+\beta_k(B+\nu_k)\notag\\
&\le
\frac{\beta_k l_g}{ \rho}
(l_f+\beta_k l_g+l_r)
+\beta_k(B+\nu_k)
\le C_k. \notag
\end{align}
\end{proof}

The following lemma establishes a separation oracle for \cref{alg:pdmcp}: whenever the residual $R(\widehat{y})$ exceeds $\epsilon_k$, the gradient $\nabla D_k(\widehat{y})$ strictly separates $\widehat{y}$ from a neighborhood of the dual optimal solution $\bar{y} := \operatorname{argmax}_{y \in \mathcal{B}_1} D_k(y)$. 

\begin{lemma}\label{lem:ctb}
Under the assumptions of \cref{lem:barbeta}, let $\epsilon_k > 0$ and $\widehat{y} \in \mathcal{B}_1$. If $R(\widehat{y}) > \epsilon_k$, then
\begin{align}\label{eq:sep_ineq}
\langle \nabla D_k(\widehat{y}), y - \widehat{y} \rangle > 0 \quad \text{for all } y \in \mathcal{B}_{\eta_k}(\bar{y}) := \{y \in \mathcal{B}_1 : \|y - \bar{y}\| \le \eta_k\},
\end{align}
where $\eta_k = \min\{\eta_k^+, 1\}$, and $\eta_k^+$ is the   positive root of
$ 
\eta + \sqrt{\frac{2\eta C_k}{\beta_k \nu_k}} = \frac{\nu_k \rho \epsilon_k}{\beta_k l_g(\nu_k \rho + \beta_k l_g^2)},
$ 
with $C_k$ defined in \cref{lem:Ck}.
\end{lemma}
 
\begin{proof}
Recall $\bar{y} = \operatorname{Proj}_{\mathcal{B}_1}\left(\frac{J_g(x^k)x(\bar{y}) + b_k}{\nu_k}\right)$ from  \cref{eq:pdopt}. Using the definition of $R(\widehat{y})$, the nonexpansiveness of $\operatorname{Proj}_{\mathcal{B}_1}$, and the bound $\|J_g(x^k)\| \le l_g$, we obtain
\begin{align*}
R(\widehat{y})
&\le \beta_k l_g \left\| \operatorname{Proj}_{\mathcal{B}_1}\left(\tfrac{J_g(x^k)x(\widehat{y}) + b_k}{\nu_k}\right) - \widehat{y} \right\| \\
&\le \beta_k l_g \left( \left\| \operatorname{Proj}_{\mathcal{B}_1}\left(\tfrac{J_g(x^k)x(\widehat{y}) + b_k}{\nu_k}\right) - \operatorname{Proj}_{\mathcal{B}_1}\left(\tfrac{J_g(x^k)x(\bar{y}) + b_k}{\nu_k}\right) \right\| + \|\widehat{y} - \bar{y}\| \right) \\
&\le \frac{\beta_k l_g^2}{\nu_k} \|x(\widehat{y}) - x(\bar{y})\| + \beta_k l_g \|\widehat{y} - \bar{y}\| \le \frac{\beta_k l_g (\nu_k \rho + \beta_k l_g^2)}{\nu_k \rho} \|\widehat{y} - \bar{y}\|,
\end{align*}
where the last inequality uses the   definition of $x(y)$ and the $1$-Lipschitz continuity of the proximal operator.
   Thus, $R(\widehat{y}) > \epsilon_k$ implies $\|\widehat{y} - \bar{y}\| > \frac{\nu_k \rho \epsilon_k}{\beta_k l_g(\nu_k \rho + \beta_k l_g^2)}$. For all $y \in \mathcal{B}_{\eta_k}(\bar{y})$, since $\eta_k \le \eta_k^+$ and $\eta_k^+$ satisfies $ 
\eta_k^+ + \sqrt{\frac{2\eta_k^+ C_k}{\beta_k \nu_k}} = \frac{\nu_k \rho \epsilon_k}{\beta_k l_g(\nu_k \rho + \beta_k l_g^2)},
$ we have
\begin{align}\label{eq:yhat-y}
\|\widehat{y} - y\| \ge  \|\widehat{y} - \bar{y}\| - \|y - \bar{y}\|  >\hspace{-1mm} \frac{\nu_k \rho \epsilon_k}{\beta_k l_g(\nu_k \rho + \beta_k l_g^2)} - \eta_k^+ \hspace{-1mm}=  \sqrt{\frac{2\eta_k^+ C_k}{\beta_k \nu_k}} \ge\hspace{-1mm} \sqrt{\frac{2\eta_k C_k}{\beta_k \nu_k}}.
\end{align}

Next, concavity of $D_k(y)$ and the Cauchy--Schwarz inequality give
\[
D_k(y) - D_k(\bar{y}) \ge \langle \nabla D_k( {y}), y - \bar{y} \rangle \ge -\|\nabla D_k( {y})\| \|y - \bar{y}\| \ge -\eta_k C_k,
\]
where the last step uses $\|\nabla D_k( {y})\| \le C_k$ from \cref{lem:Ck}. Since $D_k$ is $(\beta_k \nu_k)$-strongly concave and $D_k(\widehat{y}) \le D_k(\bar{y})$,  
combining this bound with \cref{eq:yhat-y} yields
\begin{align*}
\langle \nabla D_k(\widehat{y}), y - \widehat{y} \rangle
&\ge D_k(y) - D_k(\widehat{y}) + \frac{\beta_k \nu_k}{2} \|\widehat{y} - y\|^2 > D_k(y) - D_k(\bar{y}) + \eta_k C_k \ge 0.
\end{align*} 
\end{proof}

From this lemma, we know that if $\widehat y\in \mathcal B_{\eta_k}(\bar{ y})$, we will get $R(\widehat y)
\le  \epsilon_k$, which implies that
$\mathcal B_{\eta_k}(\bar{ y}) \subseteq \{\widehat y:R(\widehat y)
\le  \epsilon_k\}.$ Also, if $\widehat y \in \mathcal B_1$ and $R(\widehat y)>  \epsilon_k$, we can always get a separation oracle $\nabla D_k(\widehat y)$ between $\widehat y$ and $\mathcal B_{\eta_k}(\bar{ y})$. 
 
When $\widehat y$ is outside of $\mathcal B_1$, the vector $w= -\widehat y/\|\widehat y\|$ separates $\widehat y$ from $\cB_1$, satisfying $\langle w,y-\widehat y\rangle>0$ for all $y\in \mathcal B_1$. Hence cutting-plane methods can locate a point $y^t\in\mathcal B_{\eta_k}(\bar{ y})$ with $R(y^t)\le \epsilon_k$, as shown in \cref{alg:pdmcp}. We consider two classical cutting-plane updates: 
Vaidya's method~\cite{vaidya1996new,anstreicher1997vaidya} and the ellipsoid method~\cite{bland1981ellipsoid}.  
Detailed implementations 
are provided in Appendix~\ref{app:cps}.

\begin{algorithm}[h!]
\caption{Dual cutting-plane method for solving $\min_{x\in\mathbb R^d}\max_{y\in\mathcal B_1}\phi_k(x,y)$}
\label{alg:pdmcp}
\begin{algorithmic}[0]
\STATE \textbf{Input:}
$x^k\in\mathcal X$, $\beta_k>0$, $\nu_k>0$,
$\rho\ge\underline\rho$, and $\epsilon_k>0$.

\STATE Choose a cutting-plane method and set $y^0=0$.
For Vaidya's method, set \mbox{$\mathcal P^0=[-1,1]^n$.}
For the ellipsoid method, set $\mathcal P^0=\mathcal B_1$.

\FOR{$t=0,1,\ldots$}
    \IF{$\|y^t\|>1$}
        \STATE Set $w^t=-y^t/\|y^t\|$.
    \ELSE
        \STATE Compute
        $
        x(y^t)=
        \operatorname{prox}_{r/\rho+\delta_{\mathcal X}}
        (
        x^k-
        \frac{\nabla f(x^k)+\beta_kJ_g(x^k)^\top y^t}{\rho}
        ).
        $
        \STATE Compute  $R(y^t)=\beta_k \|J_g(x^k) ^\top
(\mathrm{Proj}_{\mathcal B_1}(\frac{g(x^k)+J_g(x^k)(x(y^t)-x^k)}{{\nu_k}})-y^t)\|$.
        \IF{$R(y^t)\le\epsilon_k$}
            \RETURN $x(y^t)$.
        \ENDIF
        \STATE Set
        $
        w^t=
        \beta_k\bigl(
        g(x^k)+J_g(x^k)(x(y^t)-x^k)-\nu_k y^t
        \bigr).
        $
    \ENDIF
    \STATE Apply the prescribed update of the selected
    cutting-plane method to obtain
    $(\mathcal P^{t+1},y^{t+1}\in \mathcal P^{t+1})$, maintaining
    $
    \mathcal P^t\cap
    \{y:\langle w^t,y-y^t\rangle\ge0\}
    \subseteq
    \mathcal P^{t+1}.
    $
\ENDFOR
\end{algorithmic}
\end{algorithm}

The following lemma summarizes the volume reduction rates of the two cutting-plane methods incorporated in \cref{alg:pdmcp}. 

\begin{lemma}\label{lem:volume-reduction} 
Suppose that the localization sets $\{\mathcal{P}^t\}_{t\ge 0}$ are generated by either Vaidya's method or the ellipsoid method described in Appendix~\ref{app:cps}. Then, for all $t\ge 0$, 
\begin{align}
\label{eq:volume-reduction} \frac{\operatorname{Vol}(\mathcal{P}^t)}{\operatorname{Vol}(\mathcal{P}^0)} \le \exp\bigl(\kappa n\ln(2n)-\theta t\bigr), 
\end{align} 
where $\operatorname{Vol}(\cdot)$ denotes the   $n$-dimensional volume. The parameters $\kappa$ and $\theta$ are given by $\kappa = c_1^{\mathrm V}$ and $\theta = c_2^{\mathrm V}$ for Vaidya's method, and $\kappa = 0$ and $\theta = \frac{1}{2n+2}$ for the ellipsoid method, where $c_1^{\mathrm V}$ and $c_2^{\mathrm V}>0$ are universal constants independent of $n$ and~$t$. 
\end{lemma}

 \begin{proof}
We analyze the two cutting-plane updates separately:

\paragraph{Vaidya's method} Let $m_t$ be the number of constraints defining $\cP^t$ at iteration~$t$, and $\mathbf{\omega}^t$ be the volumetric center of $\mathcal{P}^t$. Applying the  bound from \mbox{\cite[Eq.~(3.4)]{anstreicher1997vaidya}}, we have
\begin{align}
    \label{eq:volumnbound}
    \ln\tfrac{\operatorname{Vol}(\mathcal{P}^t)}{\operatorname{Vol}(\mathcal{P}^0)} \le n\ln(m_t) - V_{\mathcal{P}^t}(\mathbf{\omega}^t),
\end{align}
where $V_{\mathcal{P}^t}$ is the volumetric barrier. The  analysis in \cite[Theorems~2.6 and 3.2]{anstreicher1997vaidya} and \mbox{\cite[Corollary 6.6]{anstreicher1997vaidya}} ensures that $m_t =  {O}(n)$ and   there exist   \mbox{$c_2^{\mathrm{V}},c_3^{\mathrm{V}} > 0$} such that $V_{\mathcal{P}^t}(\mathbf{\omega}^t) \ge c^{\mathrm{V}}_2 t - c_3^{\mathrm{V}} n\ln(2n)$ for all $t \ge 0$. Substituting these estimates into \eqref{eq:volumnbound} yields
\[\ln \tfrac{\operatorname{Vol}(\mathcal{P}^t)}{\operatorname{Vol}(\mathcal{P}^0)} \le   c^{\mathrm{V}}_1 n\ln(2n) - c^{\mathrm{V}}_2 t  
\]
for a sufficiently large constant $c^{\mathrm{V}}_1 > 0$ that absorbs the $ O(n)$ term.

\paragraph{Ellipsoid method} For $n = 1$, each bisection step halves the interval length, giving $\operatorname{Vol}(\mathcal{P}^t)/\operatorname{Vol}(\mathcal{P}^0) = 2^{-t} \le \exp( -t/4 ) \le \exp( -t/(2n+2) )$. For $n \ge 2$, each central-cut update satisfies the geometric reduction rate \cite[Eq.~(2.11)]{bland1981ellipsoid}:
\[
\frac{\operatorname{Vol}(\mathcal{P}^{t+1})}{\operatorname{Vol}(\mathcal{P}^t)} = \frac{n}{n+1} \left( \frac{n^2}{n^2-1} \right)^{(n-1)/2} \le \exp\left(-\frac{1}{2(n+1)}\right).
\]
Applying this inequality   from $0$ to $t-1$ yields $\operatorname{Vol}(\mathcal{P}^t)/\operatorname{Vol}(\mathcal{P}^0) \le \exp\big(-t/(2n+2)\big)$ for all $n \ge 1$.
\end{proof}

Combining \cref{lem:ctb,lem:volume-reduction} yields an iteration bound for \cref{alg:pdmcp}.

\begin{lemma}\label{lem:cpcmp}
Under the conditions of \cref{lem:barbeta}, suppose that \cref{alg:pdmcp} employs either Vaidya's method or the ellipsoid method.
Then it terminates within
\[
T_{2}^k
:= 
\left\lceil
\frac{
\kappa n\ln(2n)
+n\ln(2\sqrt n/\eta_k)
}{\theta}
\right\rceil +1
\]
cutting-plane iterations, where $\eta_k$ is defined in \cref{lem:ctb}, and $\kappa$ and $\theta$ are specified in \cref{lem:volume-reduction}.
\end{lemma}

\begin{proof}
    Suppose that for all $t\le T_2^k$, we have $\|y^t\| >1$ or $R(y^t) > \epsilon_k$. When $\|y^t\| >1$, setting $w^t =  -y^t/\|y^t\|$ gives 
$ \langle w^t, y - y^t \rangle> 0$ for all $y\in \mathcal B_{\eta_k}(\bar y)\subseteq \mathcal B_1.
$
For the case where $y^t\in \mathcal B_1$ and $R(y^t)
>  \epsilon_k$, by setting \[w^t = \nabla D_k(y^t) = \beta_k(g(x^k)+J_g(x^k) (x(y^t)-x^k)-{\nu_k} y^t),\] we also obtain $\langle w^t, y - y^t \rangle> 0$ for all $y\in \mathcal B_{\eta_k}(\bar y)$ by \cref{lem:ctb}.
 In both cases, it holds that $\mathcal B_{\eta_k}(\bar y)\subseteq \{y:\langle w^t, y - y^t \rangle> 0\}$. Since \mbox{$
    \mathcal P^t\cap
    \{y:\langle w^t,y-y^t\rangle\ge0\}
    \subseteq
    \mathcal P^{t+1}
    $} and $\mathcal B_{\eta_k}(\bar y)\subseteq \mathcal P^0$, induction yields $\mathcal B_{\eta_k}(\bar y)\subseteq \mathcal P^t$ for all $t\le T_2^k$. Then we obtain 
    \[\begin{aligned}
     \operatorname{Vol}(\mathcal B_{\eta_k}(\bar y)) &\le \operatorname{Vol}(\mathcal P^{T_2^k}) {\le} e^{\kappa n\ln (2n)-\theta T^k_2}\operatorname{Vol}(\mathcal P^0)\\&<e^{- n\ln \frac{2\sqrt n}{\eta_k}}\operatorname{Vol}(\mathcal P^0)\le e^{- n\ln \frac{2\sqrt n}{\eta_k}} \cdot 2^n=\left(\frac{\eta_k}{\sqrt n}\right)^{n},
 \end{aligned}\]
 where the second inequality is from \cref{eq:volume-reduction}. On the other hand, since $\bar y \in \mathcal B_1$ and $\eta_k = \min\{\eta_k^+, 1\} \le 1$, the ball of radius $\eta_k/2$ centered at $(1 - \eta_k/2)\bar y$ lies in $\mathcal B_{\eta_k}(\bar y)$.
Hence,
\begin{align*}
    \operatorname{Vol}(\mathcal B_{\eta_k}(\bar y)) &= \operatorname{Vol}(\{y\in \mathcal B_1: \|y-\bar y\|\le \eta_k\}) \ge \operatorname{Vol}\left(
\left\{z:
\left\|z-\left(1-\frac{\eta_k}{2}\right)\bar y
\right\|\le\frac{\eta_k}{2}\right\}\right) \\
    &\ge   \operatorname{Vol}\left(\left\{z:\|z-\left(1-\frac{\eta_k}{2}\right)\bar y\|_{\infty}  \le \frac{\eta_k}{2\sqrt{n}}\right\}\right) = \left(\frac{\eta_k}{\sqrt{n}}\right)^{n},
\end{align*} contradicting the strict upper bound above.
\end{proof}

The following theorem establishes the total arithmetic complexity when Vaidya's method or the ellipsoid method is used within \cref{alg:spplm}.

\begin{theorem} \label{thm:dcp}
Suppose that \cref{ass:problemsetup,ass:regularity} hold. Let $\epsilon \in (0,1)$, and set $ \nu_k$ and $\epsilon_k$ as those in \cref{thmcdt} with $c_\nu\in(0,\frac{B\sqrt{2}\ln 2} 2]$ and $\epsilon_0\in(0,1]$. 
If each prox-linear subproblem in \cref{alg:spplm} is solved by \cref{alg:pdmcp}, the total arithmetic complexity to find an $\epsilon$-KKT point of problem \cref{eq:model} is:
\begin{enumerate}[leftmargin=*, labelindent=-5pt, label=(\roman*)]
    \item $\widetilde{O}\big((\mathcal{C}_{\mathrm{FO}} + dn^2 + n\mathcal{C}_r + n^{3.373})\epsilon^{-2}\big)$ when using Vaidya's method in Appendix~\ref{app:vaidya};
    \item $\widetilde{O}\big((\mathcal{C}_{\mathrm{FO}} + dn^3 + n^2\mathcal{C}_r + n^4)\epsilon^{-2}\big)$ when using the ellipsoid method in Appendix~\ref{app:ellip}.
\end{enumerate}
\end{theorem}

\begin{proof}
Under the stated assumptions, by \cref{thm:subconv}, \cref{alg:spplm} finds an $\epsilon$-KKT point within
$K=O(\epsilon^{-2})$ outer iterations.
Moreover, \cref{lem:bndcon,lem:barbeta,lm:barrho,lem:Ck} give
$\|g(x^k)\|\le B,\beta_k = \Theta(1)$, $\rho=\Theta(1)$ and $C_k=\Theta(1)$
  for all $0\le k<K$.

We first bound $\eta_k$ from below for all $0 \le k < K$. When $\eta_k^+ \le 1$, we have $\eta_k^+ \le \sqrt{\eta_k^+} = O(\sqrt{\eta_k^+/\nu_k})$ since $\nu_k = O(1)$. Combining the defining equation for $\eta_k^+$ in \cref{lem:ctb} with the parameter bounds above gives
\begin{align*} 
\Theta(\nu_k \epsilon_k) = \frac{\nu_k \rho \epsilon_k}{\beta_k l_g(\nu_k \rho + \beta_k l_g^2)} = \eta_k^+ + \sqrt{ \frac{2\eta_k^+ C_k}{\beta_k \nu_k}} = O\left(\sqrt{\frac{\eta_k^+}{\nu_k}}\right).
\end{align*}   Squaring both sides yields $\eta_k^+ = \Omega(\nu_k^3 \epsilon_k^2)$. Since $\nu_k^3 \epsilon_k^2 = O(1)$, it follows for both $\eta_k^+ \le 1$ and $\eta_k^+ > 1$ that $\eta_k = \min\{\eta_k^+, 1\} = \Omega(\nu_k^3 \epsilon_k^2)$. By~\cref{thmcdt}, 
$
    \nu_k=\widetilde\Omega(\epsilon) $ and $
    \epsilon_k=\Omega(\epsilon^2) $ for all $ 0\le k<K.
$
Consequently,
$
    \eta_k=\widetilde\Omega(\epsilon^7)
$ uniformly over these iterations. 
At each inner iteration of \cref{alg:pdmcp}, given $g(x^k)$, $\nabla f(x^k)$, and $J_g(x^k)$, evaluating $(x(y^t), R(y^t), w^t)$ requires $O(dn + \mathcal{C}_r)$ arithmetic operations. By \cref{lem:cpcmp}, we bound the per-call arithmetic cost of \cref{alg:pdmcp} at iteration $k$ as follows:
\begin{enumerate}[leftmargin=*, labelindent=-5pt, label=(\roman*)]
    \item \emph{Vaidya's method:} With $\kappa = c_1^{\mathrm{V}}$ and $\theta = c_2^{\mathrm{V}}$ in \cref{lem:volume-reduction}, \cref{alg:pdmcp} terminates in $T_2^k = O(n \ln(2n/\eta_k))   = O(n \ln(n/\epsilon)) $ inner iterations. Updating $\mathcal P^t$ and $y^t$ in \cref{alg:vaidya-update} requires $O(n^{2.373})$ operations via fast matrix multiplication~\cite{williams2024new} and inversion~\cite[Section 7.3]{anstreicher1997vaidya}. Thus, the arithmetic cost for a single invocation of \cref{alg:pdmcp}  is $\widetilde{O}((dn + \mathcal{C}_r + n^{2.373})T_2^k) = \widetilde{O}(dn^2 + n\mathcal{C}_r + n^{3.373}).$ 
    \item \emph{Ellipsoid method:} With $\kappa = 0$ and $\theta = \frac{1}{2n+2}$, \cref{alg:pdmcp} terminates in $T_2^k = O(n^2 \ln(2n/\eta_k)) = O(n^2 \ln(n/\epsilon)) $ inner iterations. Each ellipsoid update in \cref{alg:ellipsoid-dcp} costs $O(n^2)$ operations. Thus, each call to \cref{alg:pdmcp} costs \mbox{$\widetilde{O}((dn + \mathcal{C}_r + n^2)T_2^k) = \widetilde{O}(dn^3 + n^2\mathcal{C}_r + n^4)$} operations.
\end{enumerate}
Finally, multiplying the per-call subsolver cost and oracle cost $\mathcal{C}_{\mathrm{FO}}$ by the total number of invocations $(\lceil \log_\gamma(\bar{\rho}/\underline{\rho}) \rceil + 1)K = O(\epsilon^{-2})$ yields both bounds.
\end{proof} 

Vaidya's method has a better polynomial dependence on $n$ than the ellipsoid method.
We therefore compare Vaidya's method with Nesterov's APG subsolver in \cref{thm:cmpapg}.
Both methods share the first-order oracle cost $O(\mathcal C_{\mathrm{FO}}\epsilon^{-2})$, while the APG and Vaidya subsolvers incur additional arithmetic costs of
\[
\widetilde  O\left((dn+\mathcal C_r)\epsilon^{-5/2}\right)
\quad\text{and}\quad
 \widetilde  O\left((dn^2+n\mathcal C_r+n^{3.373})\epsilon^{-2}\right),
\]
respectively. 
Suppressing the polylogarithmic factors, Vaidya’s method has lower subproblem cost whenever
 \mbox{$ n+ {n^{3.373}}/(dn+\mathcal C_r)=o(\epsilon^{-1/2}) $}.  
 With the other quantities fixed, this condition becomes easier to satisfy as $n$ decreases, $d$ increases, or $\epsilon$ decreases. Conversely, when
 $\epsilon^{-1/2}=o(n)$,
 Nesterov's APG subsolver has lower arithmetic cost.

 More advanced cutting-plane methods \cite{lee2015faster,jiang2020improved} can improve the polynomial dependence $\widetilde{O}(n^{3.373})$ via randomized leverage-score maintenance.
We do not pursue such implementations, since they require substantially more involved algorithmic machinery than the volumetric and ellipsoid routines adopted here.

\section{Numerical experiments}
\label{sec:numerics}

In this section, we evaluate two variants of \cref{alg:spplm}: SEPi-APG  (\cref{alg:acceleratedNesterov}) and SEPi-DCP (\cref{alg:pdmcp}).
We compare both variants with 
iPPP~\cite{lin2022complexity}
and iALM~\cite{sahin2019inexact}.   
All experiments are run in \mbox{MATLAB R2025b} on an  Apple M4 MacBook Air laptop with 24 GB RAM.

 Throughout this section, the data pass count measures the number of full data evaluations. 
For SEPi-APG and SEPi-DCP, the line search for each  trial $\rho$ counts for  
one data pass.  For iPPP and iALM, data passes are accumulated from inner function and gradient evaluations. For all methods, stationarity violations are measured by the same KKT residual of the original constrained problem.

\subsection{Tested problem}  
\label{subsec:mnpc-numerics}

We 
consider the multi-class Neyman--Pearson classification (mNPC) 
in  
\cite{lin2022complexity}, using the datasets \texttt{MNIST}~\cite{lecun1998mnist} and \texttt{covtype}~\cite{covertype_31}.  
 For both datasets, all features are scaled to $[0,1]$, and each sample vector is appended with a constant~$1$.

Let $K$ denote the number of classes in this section and let $D_k=\{a_i^{(k)}\}_{i=1}^{n_k} \subseteq\RR^p$ denote the samples in class $k$.  We write $x=[x_1,\ldots,x_K]\in\RR^{p\times K}$ and impose the block constraint
$\mathcal X=\{x:\ \|x_k\|\le 0.3,\ k=1,\ldots,K\}$. 
For class $k$, define 
\[
    \ell_k(x)
    =
    \frac{1}{n_k}\sum_{j\ne k}\sum_{i=1}^{n_k}
    \frac{1}{1+\exp(\langle a_i^{(k)},x_k\rangle-\langle a_i^{(k)},x_j\rangle)} .
\]
The mNPC problem minimizes the  classification  loss of a target class while controlling the losses of all other classes within prescribed thresholds. 
Taking class 1 as the target class, we solve
\begin{align}
\label{eq:mnpc-exp}
    \min_{x\in\mathcal X}\ \ell_1(x),
    \quad
    \st \ \ell_k(x) \le {(K-1)}/{2},\quad k=2,\ldots,K.
\end{align}
The \texttt{MNIST} and \texttt{covtype} datasets respectively 
have $(p,K,\sum_{k=1}^K n_k)=(781,10,60000)$ and $(55,7,581012)$. 
Introducing slack variables $s \in \RR_+^{K-1}$ yields the reformulation
\begin{align}\label{newcon} 
    q(x)+s=0,\quad (x,s) \in \mathcal X \times \RR^{K-1}_+,
\end{align}
with $
    q(x):= [ \ell_2(x)-{(K-1)}/{2}; \cdots ; \ell_K(x)-{(K-1)}/{2}]$.
 For a simplified setting as discussed in \mbox{\cite[Appendix~2]{lin2022complexity}},  the constraint system in \cref{eq:mnpc-exp} satisfies \mbox{\cite[Assumption~4]{lin2022complexity}} under mild data conditions.
Whenever \cite[Assumption~4]{lin2022complexity} holds, the slack reformulation \eqref{newcon} satisfies \cref{ass:regularity}; see~\cite[Section~3.1]{zuo2026cq}.

 We  measure feasibility violation by  $\operatorname{feas}(x):=\|[q(x)]_+\|$ and stationarity residual by 
\begin{align}
\label{eq:mnpc-stat-exp}
    \operatorname{stat}(x)
    :=
    \min_{\lambda\in\RR_+^{|I(x)|}}
    \textstyle \dist\left(\nabla \ell_1(x)+\sum_{i\in I(x)}\lambda_i\nabla q_i(x),-\mathcal N_{\mathcal X}(x)\right),
\end{align}
where $I(x):=\{i:\ q_i(x)\ge0\}$. 

For SEPi-APG, we use $x^k$ as the initial point of Nesterov's APG method for solving the $k$-th subproblem.  
For SEPi-DCP, we employ the ellipsoid update \cite{bland1981ellipsoid} (\cref{alg:ellipsoid-dcp}) with the initial point $y^0 = 0$ as specified in \cref{alg:pdmcp}.   
With $n = O(1)$ equality constraints in \cref{newcon}, the ellipsoid method has the same asymptotic complexity as Vaidya's method in \cref{thm:dcp}.
Its bound in \cref{lem:volume-reduction} is free of the $\kappa n\ln(2n)$ term ($\kappa = 0$) and has $\theta = 1/(2n+2)$.
By contrast, Vaidya's method has a small constant $\theta = c_2^{\mathrm{V}}$~\cite[Corollary~6.6 and Section~7.1]{anstreicher1997vaidya}, giving a large factor~$1/\theta$ in \cref{lem:cpcmp}.
At each outer iteration, both SEPi-Prox variants perform the 
line search in \cref{alg:spplm} with $\underline{\rho}=10$ and $\gamma= 1.5$. We set $\epsilon_0=0.001$ in~\cref{thmcdt} and choose $L_0=1$ and $\gamma_{\rm u}=\gamma_{\rm d}=2$ in~\cref{alg:acceleratedNesterov}.
We fix $\beta_k\equiv\beta$ and treat $\beta$ and $c_\nu$ in~\cref{thmcdt} as tunable parameters. By~\cref{sfclgb}, setting $\beta_k\equiv\beta$ is consistent with~\cref{alg:spplm} whenever $\beta$ is large enough.

\subsection{Parameter sensitivity analysis}
\label{subsec:num-sensitivity}

In this subsection, we test SEPi-APG and SEPi-DCP  with different choices of 
the penalty parameter $\beta$ and the smoothing coefficient $c_\nu$.  
We first vary  
\mbox{$\beta\in\{1,10,100,1000\}$} with $c_\nu=0.1$.
We next vary 
\mbox{$c_\nu\in\{0.01,0.1,1,10\}$} with $\beta=10$ on
\texttt{MNIST} and $\beta=1000$ on \texttt{covtype}. 
SEPi-APG and SEPi-DCP use the same parameters on each dataset.

\begin{figure}[!htbp]
\centering 
\includegraphics[width=.24\linewidth]{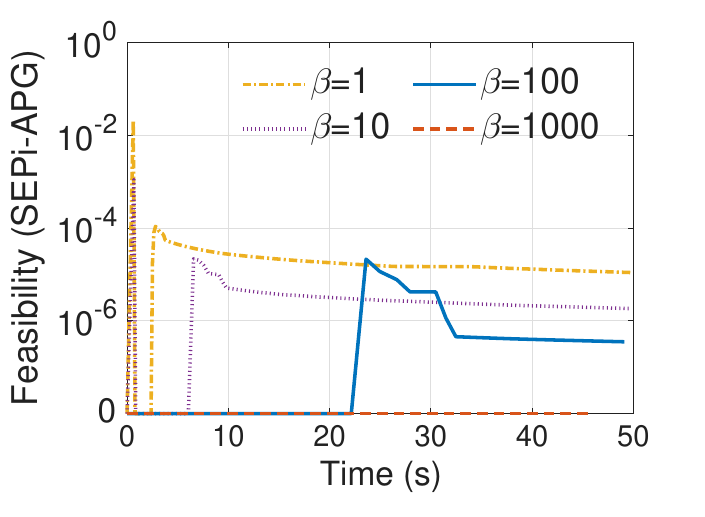}\hspace{-.008\linewidth}%
\includegraphics[width=.24\linewidth]{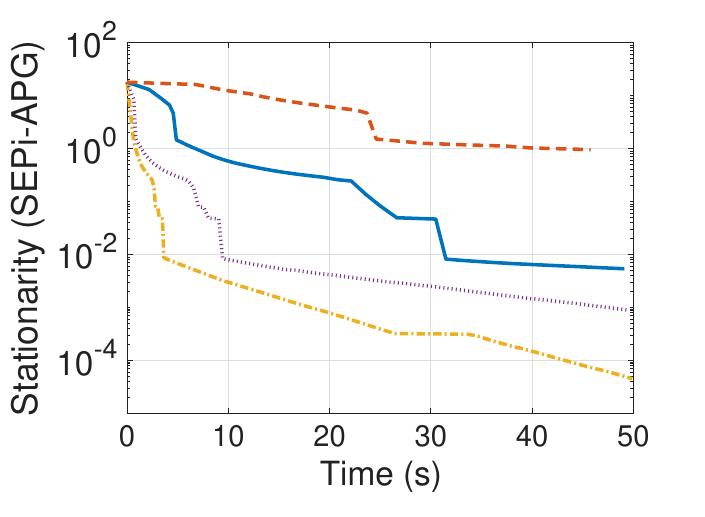}\hspace{-.008\linewidth}%
\includegraphics[width=.24\linewidth]{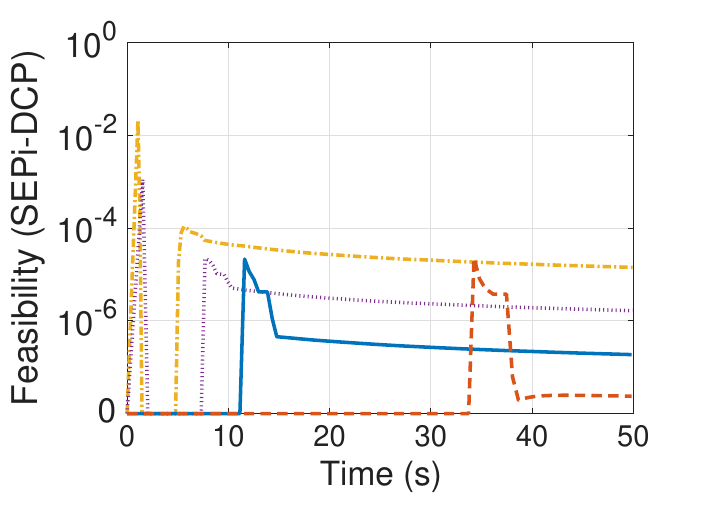}\hspace{-.008\linewidth}%
\includegraphics[width=.24\linewidth]{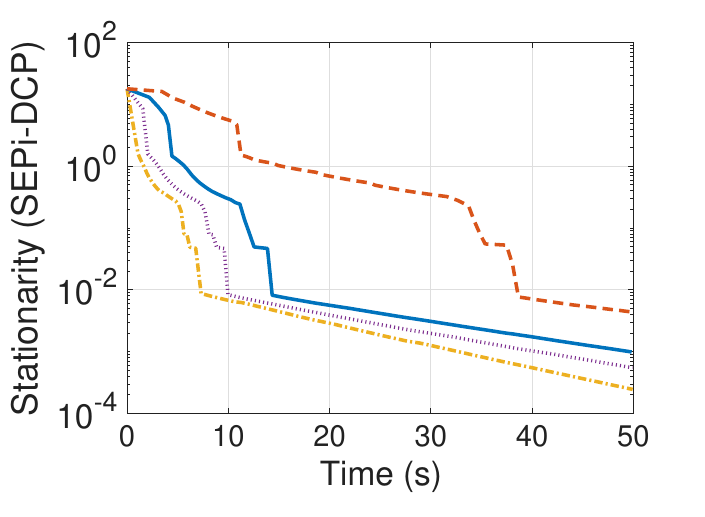}\\[-.6ex]
\hspace{.003\linewidth}
\includegraphics[width=.24\linewidth]{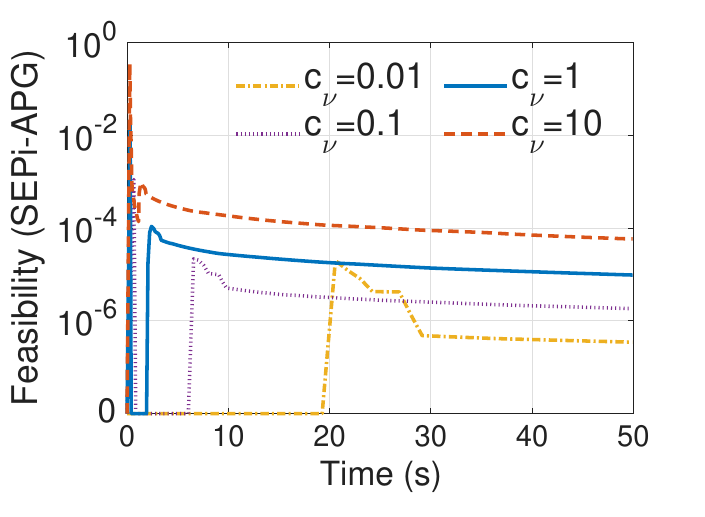}\hspace{-.008\linewidth}%
\includegraphics[width=.24\linewidth]{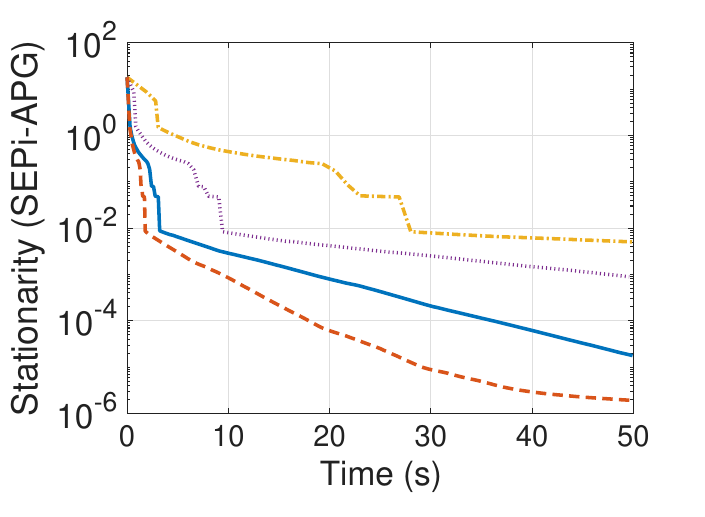}\hspace{-.008\linewidth}%
\includegraphics[width=.24\linewidth]{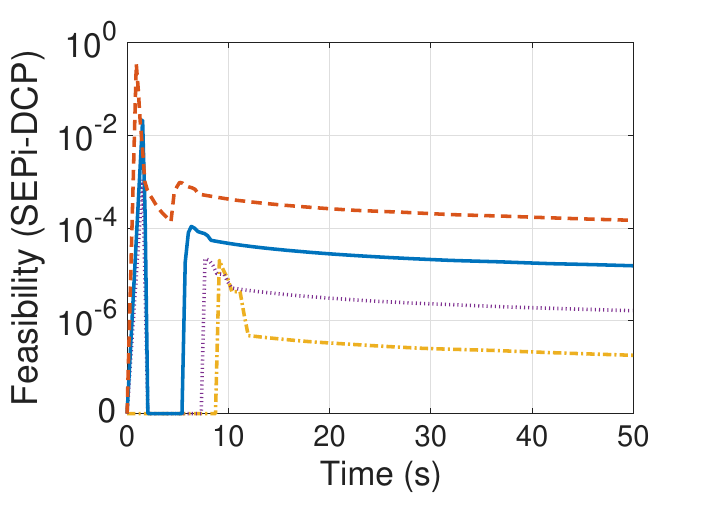}\hspace{-.008\linewidth}%
\includegraphics[width=.24\linewidth]{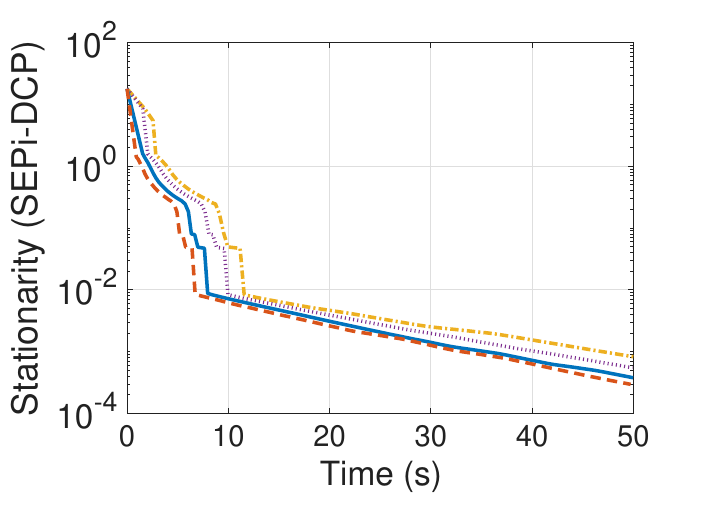}
\vspace{-3mm}
\caption{Parameter sensitivity on \texttt{MNIST}.
Rows: varying $\beta$ while fixing $c_\nu=0.1$, and varying $c_\nu$ while fixing $\beta=10$. Columns: feasibility and stationarity for SEPi-APG and SEPi-DCP.}
\label{fig:num-sens-mnist}
\vspace{-4mm}
\end{figure}

\begin{figure}[!htbp]
\centering
\includegraphics[width=.24\linewidth]{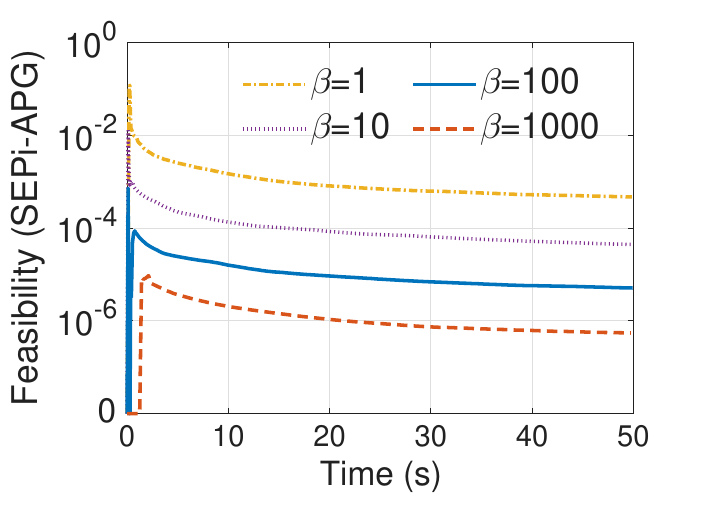}\hspace{-.008\linewidth}%
\includegraphics[width=.24\linewidth]{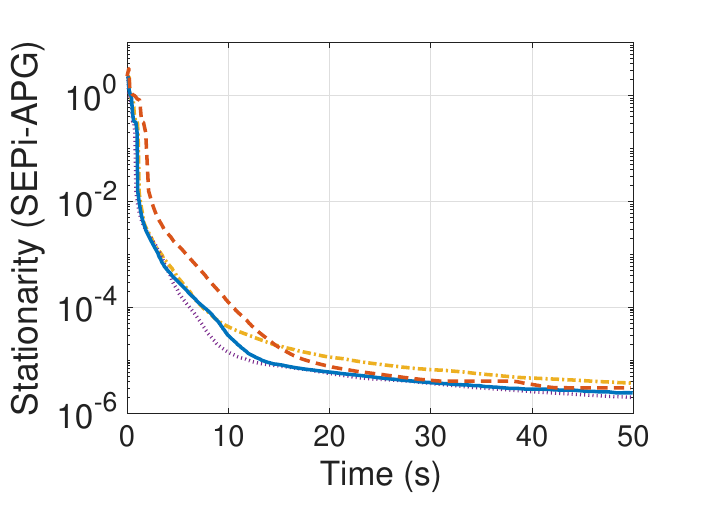}\hspace{-.008\linewidth}%
\includegraphics[width=.24\linewidth]{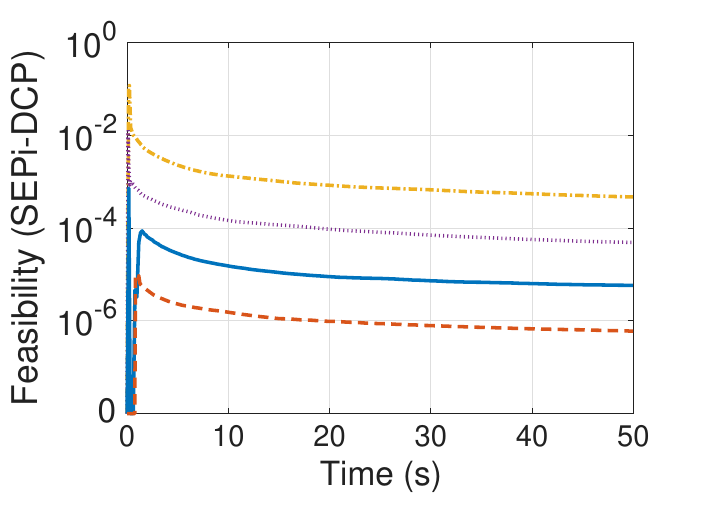}\hspace{-.008\linewidth}%
\includegraphics[width=.24\linewidth]{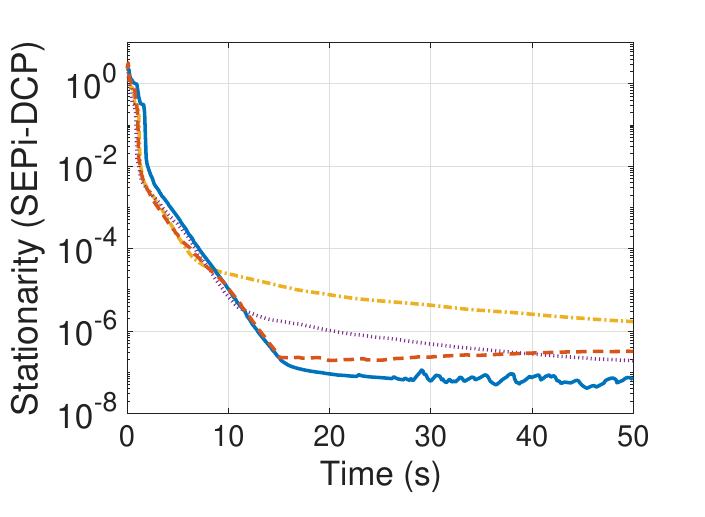}\\[-.6ex]
\hspace{.003\linewidth}
\includegraphics[width=.24\linewidth]{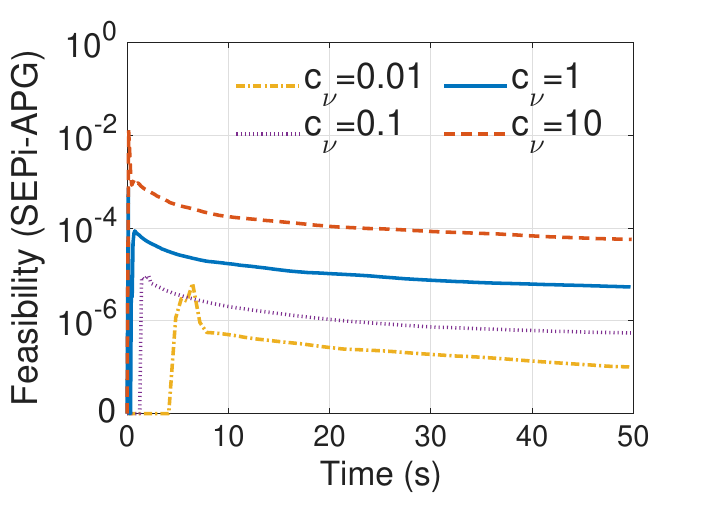}\hspace{-.008\linewidth}%
\includegraphics[width=.24\linewidth]{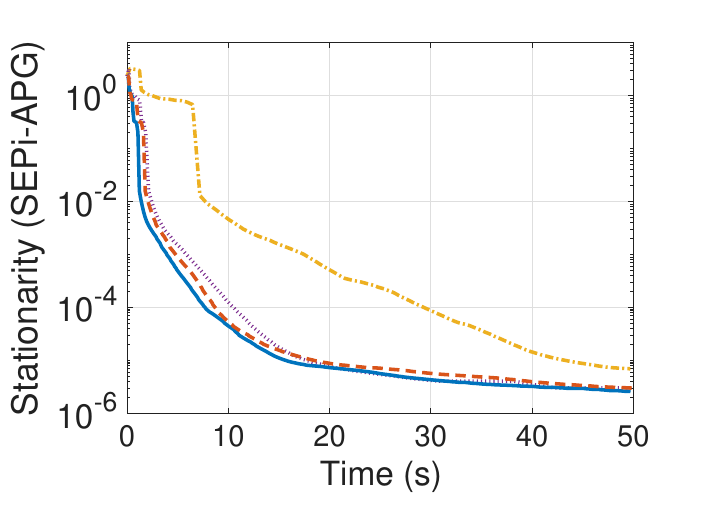}\hspace{-.008\linewidth}%
\includegraphics[width=.24\linewidth]{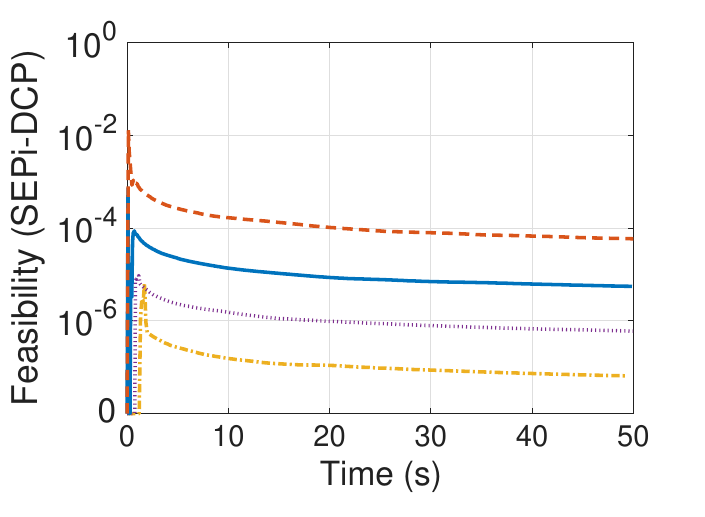}\hspace{-.008\linewidth}%
\includegraphics[width=.24\linewidth]{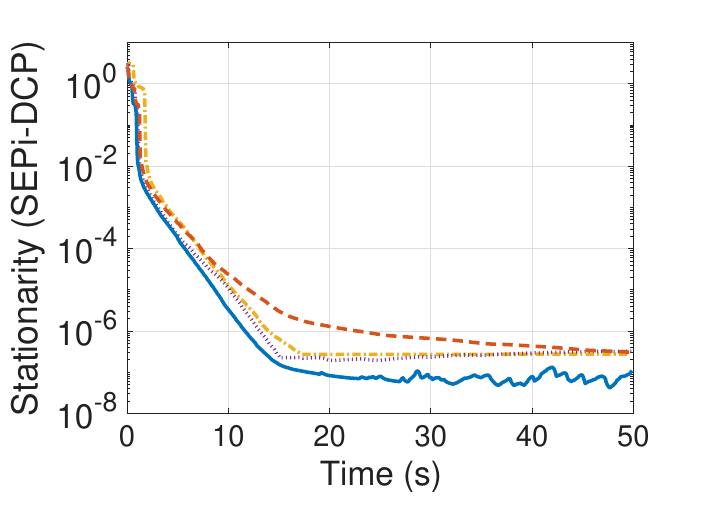}
\vspace{-3mm}
\caption{Parameter sensitivity on \texttt{covtype}.
Rows: varying $\beta$ while fixing $c_\nu=0.1$, and varying $c_\nu$ while fixing $\beta=1000$. Columns: feasibility and stationarity for SEPi-APG and SEPi-DCP.}
\label{fig:num-sens-covtype}
\end{figure}

We plot the iterates against elapsed time, as shown in Figures~\ref{fig:num-sens-mnist}
and~\ref{fig:num-sens-covtype}. From top to bottom, the rows correspond to varying $\beta$ and $c_\nu$, respectively. The first
two and last
two columns show feasibility and stationarity for SEPi-APG and SEPi-DCP, respectively. Figures~\ref{fig:num-sens-mnist} and~\ref{fig:num-sens-covtype} illustrate
the trade-offs between feasibility and stationarity as the parameters
vary. While decreasing $c_\nu$ consistently improves feasibility, an overly small value (such as $c_\nu=0.01$) degrades stationarity convergence. We therefore fix $c_\nu=0.1$ for both variants across the two datasets to balance these criteria.

The penalty parameter $\beta$ has a pronounced effect, and
its influence differs between the two problems. Increasing $\beta$
  reduces feasibility violations, but an excessively large value slows stationarity decay on \texttt{MNIST}. On \texttt{covtype}, larger $\beta$ improves constraint satisfaction while maintaining fast stationarity decay. Based on these observations we set $\beta=10$ on \texttt{MNIST} and $\beta=1000$ on \texttt{covtype} in the subsequent experiments.


\subsection{Comparisons with two baseline methods}

\begin{figure}[!htbp]
    \centering
    \includegraphics[width=0.28\textwidth]{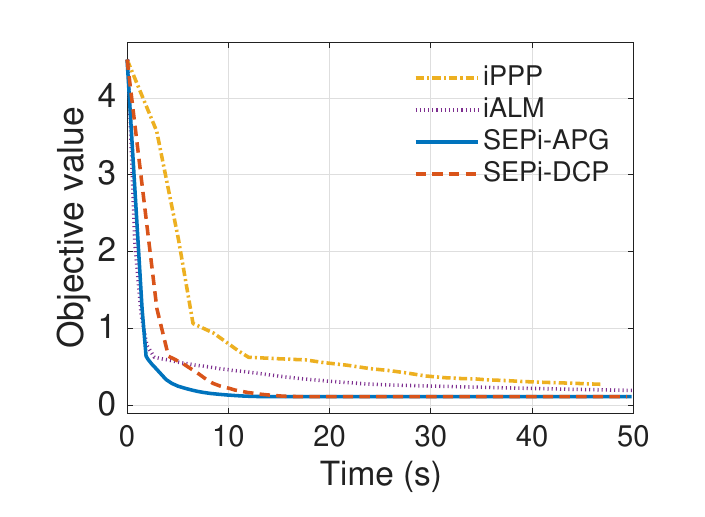} 
    \hspace{-3mm}\includegraphics[width=0.28\textwidth,viewport=0 0 338 248,clip]{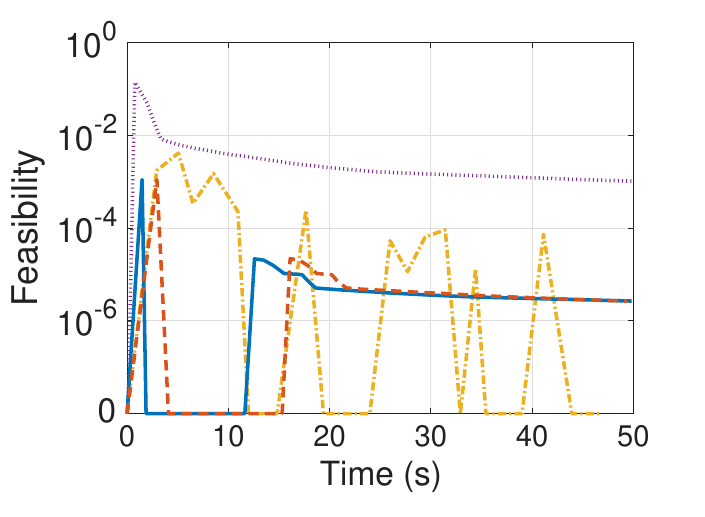} 
    \hspace{-3mm}\includegraphics[width=0.28\textwidth,viewport=0 0 338 248,clip]{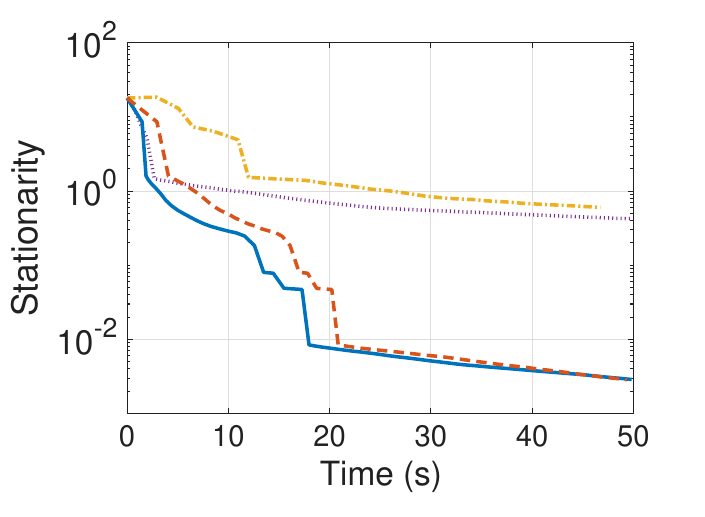}
    \vspace{-1.5mm}\\
    \includegraphics[width=0.28\textwidth,viewport=0 0 338 248,clip]{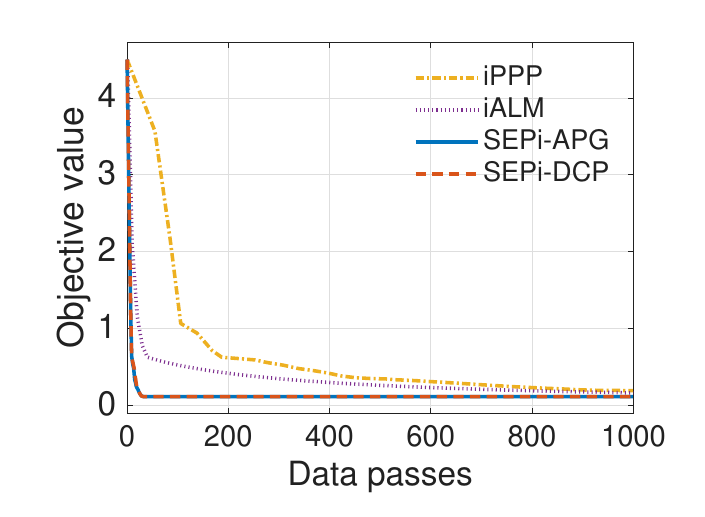} 
    \hspace{-3mm}\includegraphics[width=0.28\textwidth,viewport=0 0 338 248,clip]{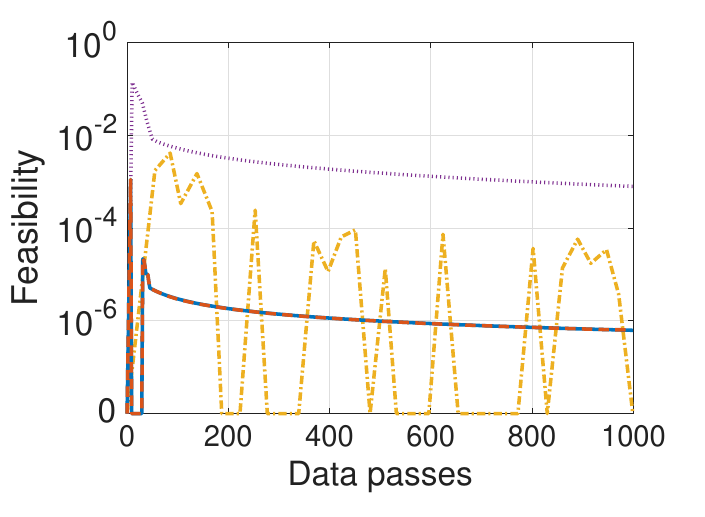} 
    \hspace{-3mm}\includegraphics[width=0.28\textwidth,viewport=0 0 338 248,clip]{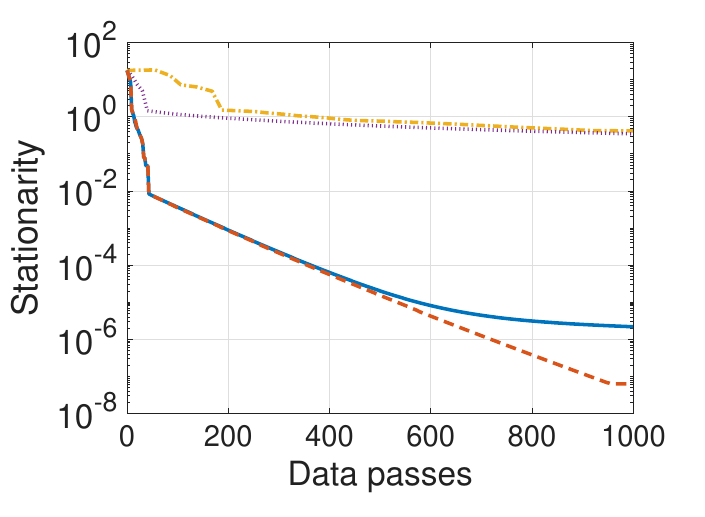}
\vspace{-3mm}
    \caption{Objective value, feasibility violation, and the stationarity residual along wall-clock time and data passes by iPPP, iALM, SEPi-APG, and SEPi-DCP for solving mNPC on the \texttt{MNIST} dataset.}
    \label{fig:mnpc-mnist-trajectories}
\vspace{-4mm}
\end{figure}

  For both SEPi-APG and SEPi-DCP, we set $c_\nu = 0.1$, with $\beta = 10$ on \texttt{MNIST} and $\beta = 1000$ on \texttt{covtype}. The baselines iPPP (with an increasing penalty parameter) and iALM are run with the same data pass budget and with the parameters set as in \cite{lin2022complexity}. All methods start from $x^0 = 0$ and run under common budgets of 50 seconds and 1000 data passes.  Comparisons of the objective values, feasibility violation, and stationarity residual on the \texttt{MNIST} dataset are plotted in \Cref{fig:mnpc-mnist-trajectories}. 
In terms of wall-clock time, SEPi-APG decays fastest in objective and stationarity residual while keeping feasibility violation small.   
In terms of data passes, both SEPi-APG and SEPi-DCP outperform the baselines, with SEPi-DCP achieving the lowest stationarity residual.

\begin{figure}[!htbp]
    \centering 
    \includegraphics[width=0.28\textwidth,viewport=0 0 338 248,clip]{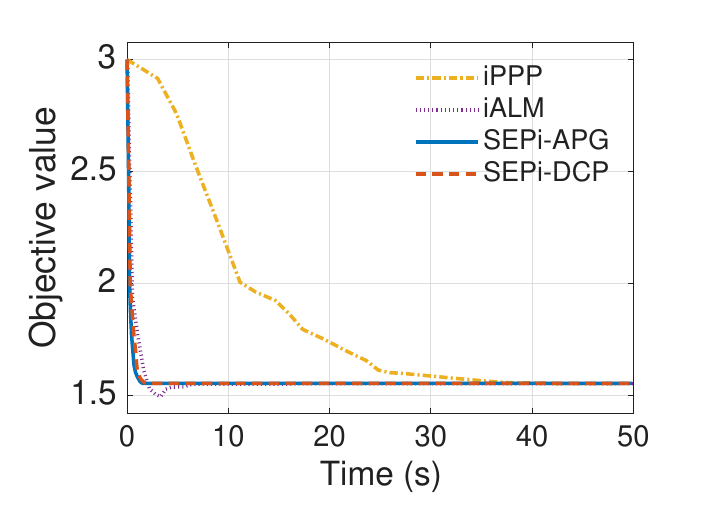} 
    \hspace{-3mm}\includegraphics[width=0.28\textwidth,viewport=0 0 338 248,clip]{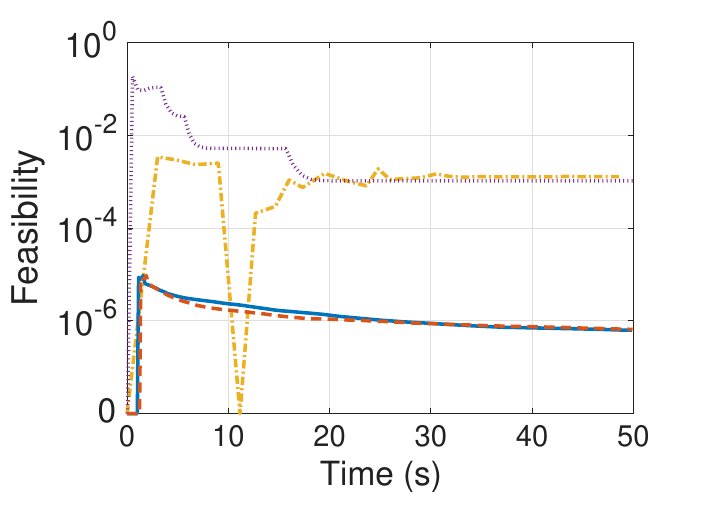} 
    \hspace{-3mm}\includegraphics[width=0.28\textwidth,viewport=0 0 338 248,clip]{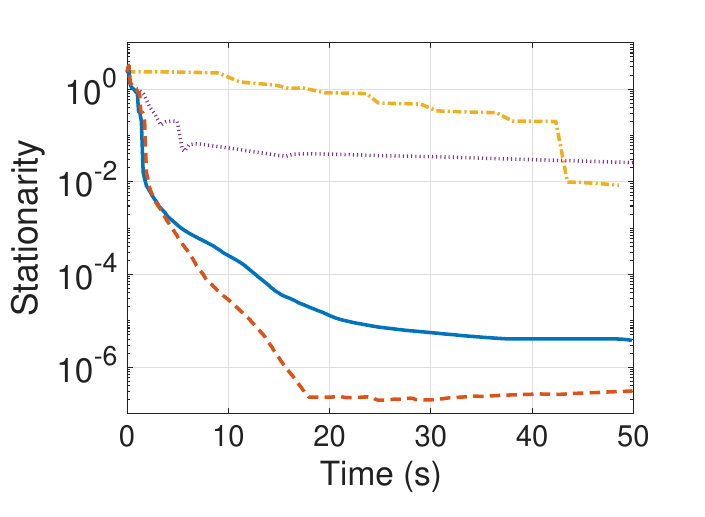}
    \vspace{-1.5mm}\\
    \includegraphics[width=0.28\textwidth,viewport=0 0 338 248,clip]{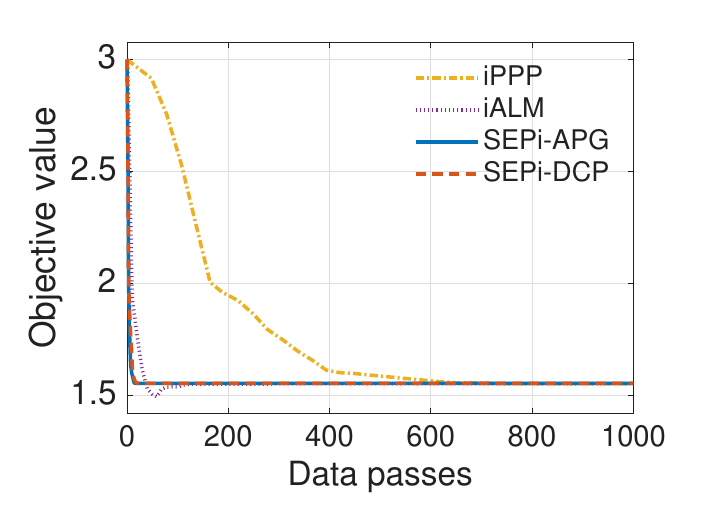} 
    \hspace{-3mm}\includegraphics[width=0.28\textwidth,viewport=0 0 338 248,clip]{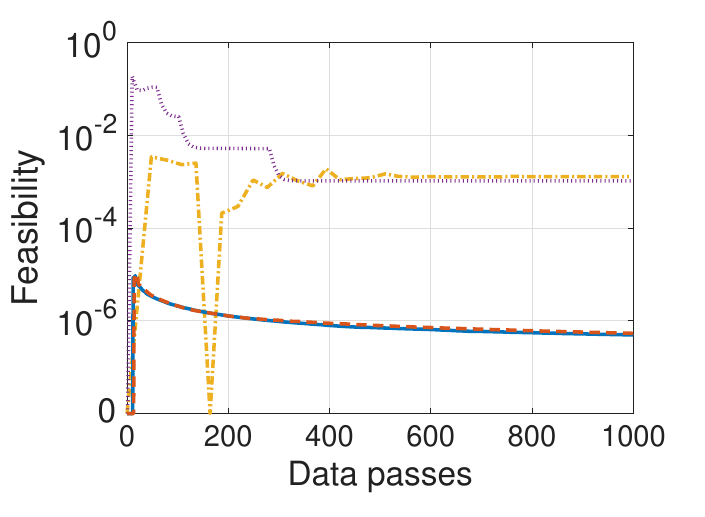} 
    \hspace{-3mm}\includegraphics[width=0.28\textwidth,viewport=0 0 338 248,clip]{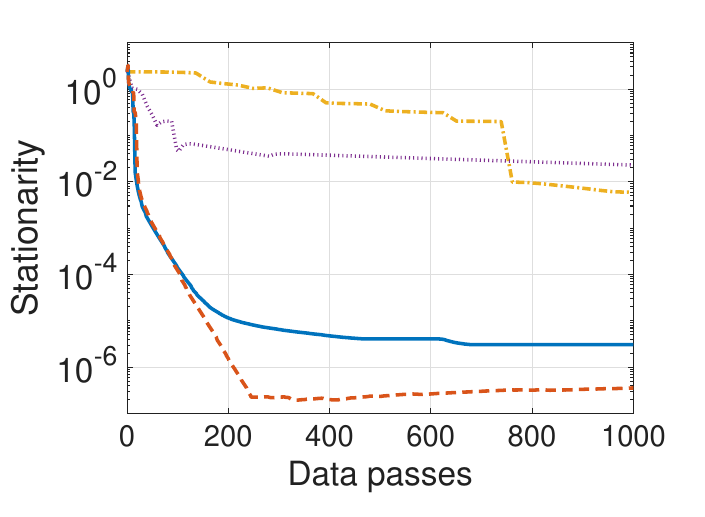}
\vspace{-3mm}
    \caption{Objective value, feasibility violation, and the stationarity residual along wall-clock time and data passes by iPPP, iALM, SEPi-APG, and SEPi-DCP for solving mNPC on the \texttt{covtype} dataset.}
    \label{fig:mnpc-covtype-trajectories}
\end{figure}

Results on the \texttt{covtype} dataset are shown in \Cref{fig:mnpc-covtype-trajectories}.  
Both SEPi-APG and SEPi-DCP outperform iPPP and iALM in terms of both wall-clock time and data passes. The two variants of SEPi-Prox reach comparable objective values and feasibility violations, while SEPi-DCP attains the smallest stationarity residual. 


\FloatBarrier

\appendix

\section{Some auxiliary results}\label{apdaux}

\begin{proof}[Proof of \cref{eq:smoothing-variation}] By the definition of $h_\nu$ in \eqref{dualhv}, we have $h_{\nu_2}(z)\geq h_{\nu_1}(z)$ for all $0<\nu_2 \le \nu_1$. Suppose that $h_{\nu_2}(z) = \langle y_2,z\rangle-\frac{\nu_2}{2}\|y_2\|^2$ for some $y_2\in\cB_1$. 
    Then, 
    \begin{align*}
        h_{\nu_2}(z)-h_{\nu_1}(z)
        &\leq
        \left(\langle y_2,z\rangle-\frac{\nu_2}{2}\|y_2\|^2\right)
        -\left(\langle y_2,z\rangle-\frac{\nu_1}{2}\|y_2\|^2\right)\\
        &=\frac{\nu_1-\nu_2}{2}\|y_2\|^2
        \leq \frac{1}{2}(\nu_1-\nu_2).
    \end{align*}
    \end{proof}

    \begin{lemma}\label{lem:dxnu}
     Consider \(\varphi(u)=\frac{u}{\sqrt{u+2}\ln(u+2)}\). It holds that \(\varphi(u)\ge a\) for all  \(a>0\) and
     $u\ge 6+32a^2\ln^2\!\bigl(4a+1\bigr)$.
    \end{lemma}

    \begin{proof}
    It is straightforward to verify that for \(u\ge6\),
    \[
    \frac{\varphi'(u)}{\varphi(u)}
    =
    \frac1 u-\frac1{2(u+2)}
    -\frac1{(u+2)\ln(u+2)}
    >
    \frac1 u-\frac1{u+2}>0,
    \]
    where we have used \(\ln(u+2)\ge\ln8>2\). Hence, \(\varphi\) is strictly increasing on \([6,\infty)\).
    
    When \(0<a\le1\) and  
    \(u\ge6\), we have
    \(\varphi(u)\ge \varphi(6)=6/(\sqrt8\ln8)>1\ge a\).
Now suppose   $a\ge1$. Set
$
A=32a^2\ln^2(4a).
$
Then
$A>6,$ and $
u\ge6+32a^2\ln^2\!\bigl(4a+1\bigr)>A.
$ Also
we have \(A+2\le2A\), and hence
    \(
    \varphi(A)\ge \frac{\sqrt A}{\sqrt2\,\ln(2A)}.
    \)
    Moreover,
    \[
    \ln(2A)
    =\ln\!\left(64a^2\ln^2(4a)\right)=2\ln(4a)+\ln4+2\ln\ln(4a)\le4\ln(4a),
    \]
    where the last inequality follows from
    \(\ln t\le t-1\) for \(t>0\) and \(\ln4<2\).
    Therefore,
    \(
    \varphi(u) > \varphi(A)\ge
    \frac{4\sqrt2\,a\ln(4a)}
    {4\sqrt2\,\ln(4a)}
    =a.
    \)  
    \end{proof}

\section{Accelerated proximal-gradient method}\label{apdx:apg}

To solve \eqref{eq:sub}, we use Nesterov's APG method \cite{nesterov2013gradient} for  problems of the form
\begin{align}\label{pAPG}
\min_{x\in \mathbb{R}^d} \ \widehat\phi(x) := \widehat f(x) + \widehat \Psi(x),
\end{align} 
where \(\widehat f\) is convex with a Lipschitz continuous gradient, and \(\widehat\Psi(x)\) is proper closed and
\(\mu\)-strongly convex. 
The method in~\cref{alg:acceleratedNesterov} adapts Nesterov's APG method with backtracking line
search~\cite[Algorithm~(4.9)]{nesterov2013gradient}.
We terminate the method when the residual $\|\nabla \widehat f(x^{t+1})-\nabla \widehat f(y^t) - M_t(x^{t+1}-y^t)\|$ reaches the prescribed tolerance $\epsilon>0$.

\begin{algorithm}[htbp!]
			\caption{{Nesterov's APG method}~\cite{nesterov2013gradient} for problem \eqref{pAPG}} 
            \label{alg:acceleratedNesterov}
            \begin{algorithmic}[0]
				\STATE \textbf{Initialization:} Choose $x^0\in\operatorname{dom}\widehat\Psi$,
$L_0>0$, $\gamma_d\ge 1$, $\gamma_u>1$, and $\epsilon>0$.
Set $L=L_0$,
$\psi_0(x)=\frac12\|x-x^0\|^2$,
$A_0=0$, and $v^0=x^0$.
				
				\FOR{$t=0,1,\ldots $}

                    \WHILE{true} 
                    \STATE  
						Choose $a>0$ satisfying $\frac{a^2}{A_t + a} = 2 \frac{1 + \mu A_t}{L}$, and set  $y = \frac{A_tx^t + av^t}{A_t + a}$. 
                        \STATE Set $z=\operatorname{prox}_{\widehat\Psi/L}(y-\frac{1}{L}\nabla\widehat f(y))$ and $\xi = \nabla \widehat{f}(z) - \nabla \widehat{f}(y) - L(z-y)$.
                    \STATE \textbf{if} $\left\langle \xi, y - z \right\rangle \ge  \frac{1}{{L}} \|\xi\|^2$,
                    \textbf{then  break}. 
                    \STATE Set $L \gets {L}\gamma_u$.
                    \ENDWHILE
                
					\STATE Set  $ y^t = y, x^{t+1} = z$, $M_t = L$, $L = \frac{M_t}{\gamma_d}$,  and
					$A_{t+1} = A_t + a$. 
                    
                    \STATE Define $\psi_{t+1}(x) := \psi_{t}(x) + a [ \widehat{f}(x^{t+1}) + \langle \nabla \widehat{f}(x^{t+1}), x - x^{t+1} \rangle + \widehat\Psi(x)]$ and set $ v^{t+1} = \argmin_x \psi_{t+1}(x)$. 
                    \STATE \textbf{if} $\|\nabla \widehat f(x^{t+1})-\nabla \widehat f(y^t) - M_t(x^{t+1}-y^t)\|\leq \epsilon$ \textbf{then} output $x^{t+1}$ and stop.
			\ENDFOR
\end{algorithmic}
		\end{algorithm}

\begin{proof}[Proof of \cref{lem:T}]

    Let $x^{k,*}$ be the optimal solution of the strongly convex subproblem \cref{eq:sub}. 
    We apply~\cref{alg:acceleratedNesterov} with
$x^0=x^k$, $\mu=\rho$, and $0<L_0\le\widetilde L_k$.
  Recall that
$
\nabla m^k(x)
=
\nabla f(x^k)
+
\beta_k J_g(x^k)^\top
\operatorname{Proj}_{\mathcal B_1}
\left(
\frac{g(x^k)+J_g(x^k)(x-x^k)}{\nu_k}
\right).
$
Thus,
$\nabla\widehat f=\nabla m^k$ has a Lipschitz constant bounded by
$\beta_k\|J_g(x^k)\|^2/\nu_k\le\widetilde L_k$.

    Let $x^t, v^t, y^t, M_t, A_t$ and $\psi_t$ be 
    generated by \cref{alg:acceleratedNesterov}, and $\psi_t^*=\min_x\psi_t(x)$.
Since $\psi_t$ is $(1+\rho A_t)$-strongly convex,
for $t\ge 1$ we have
    \begin{align*}
        &\|x^{k,*}-v^t\| \leq  \sqrt{\frac{2(\psi_t(x^{k,*})-\psi^*_t)}{\rho A_t +1}} \le \sqrt{\frac{2(A_t \widehat\phi(x^t) -\psi^*_t)+\|x^{k,*}-x^k\|^2}{\rho A_t +1}}\\
        &\le \frac{\|x^{k,*}-x^k\|}{ \sqrt{ \rho A_t +1}}\le \frac{\|x^{k,*}-x^k\|}{ \sqrt{ \rho A_t}} \le \sqrt{\frac{\gamma_u \widetilde L_k}{2\rho}} \|x^{k}-x^{k,*}\|\left[1+\sqrt{\frac{\rho}{2\gamma_u \widetilde L_k}}\right]^{-t+1},
    \end{align*}
    where the second and third inequalities use the relations in~\cite[Eq.~(4.4)]{nesterov2013gradient}
and the optimality of $x^{k,*}$, while the last inequality follows
from the growth bound for $A_t$
\mbox{in~\cite[Eq.~(4.16)]{nesterov2013gradient}}.
From the $\rho$-strong convexity of $\widehat \phi(\cdot)=F^k(\cdot;\rho)$ and \mbox{\cite[Theorem 6]{nesterov2013gradient}}, we obtain
\[\|x^{t+1}-x^{k,*}\|\le \sqrt{\frac{2(\widehat\phi(x^{t+1})-\widehat\phi(x^{k,*}))}{\rho}} \le\sqrt{\frac{\gamma_u \widetilde L_k}{2\rho}} \|x^{k}-x^{k,*}\|\left[1+\sqrt{\frac{\rho}{2\gamma_u \widetilde L_k}}\right]^{-t}.\]
Since $y^t$ is a convex combination of $v^t$ and $x^t$,
for $t\ge 1$ we obtain
\[\|y^t-x^{k,*}\| \le \max\{\|v^t-x^{k,*}\|,\|x^{t}-x^{k,*}\|\}\le \sqrt{\frac{\gamma_u \widetilde L_k}{2\rho}} \|x^{k}-x^{k,*}\|\left[1+\sqrt{\frac{\rho}{2\gamma_u \widetilde L_k}}\right]^{-t+1}.\]
    Combining these estimates with
$M_t\le\gamma_u\widetilde L_k$
from~\cite[Eq.~(4.11)]{nesterov2013gradient}
yields, for $t\ge 1$,
    \begin{equation}\label{eq:bndxk-xk*}
    \begin{aligned}
    &\|\nabla \widehat f(x^{t+1})-\nabla \widehat f(y^t) - M_t(x^{t+1}-y^t)\| \le (\widetilde L
    _k+M_t)\|x^{t+1}-y^t\|\\
    &\le (1+\gamma_u)\widetilde L_k\|x^{t+1}-y^t\| \leq (1+\gamma_u)\widetilde L_k(\|x^{t+1}-x^{k,*}\| +\|y^t-x^{k,*}\|)\\
        &\leq (1+\gamma_u)\widetilde L_k \sqrt{\frac{2\gamma_u \widetilde L_k}{\rho}}  \left\|x^{k}-x^{k,*}\right\|\left[1+\sqrt{\frac{\rho}{2\gamma_u \widetilde L_k}}\right]^{-t+1}.
    \end{aligned}\end{equation}
 The bound in~\eqref{eq:bndxk-xk*} also holds for $t=0$,
as follows from $y^0=x^k$ and the preceding estimate for $x^1$.

 Fix any $\xi^k\in\partial r(x^k)$ and let $\zeta = \nabla f(x^k) + \beta_k J_g(x^k)^\top \mathrm{Proj}_{\mathcal B_1}(g(x^k)/\nu_k) + \xi^k$. Since $0\in\mathcal{N}_{\cX}(x^k)$, we have $\zeta\in\partial F^k(x^k;\rho)$. Moreover, as \mbox{$0\in\partial F^k(x^{k,*};\rho)$}, the \mbox{$\rho$-strong} monotonicity of $\partial F^k(\cdot;\rho)$ yields
\mbox{$ \rho\|x^k-x^{k,*}\|^2
    \le
    \langle\zeta,x^k-x^{k,*}\rangle
    \le
    \|\zeta\|\|x^k-x^{k,*}\|. $}
Consequently,
\begin{align*}
    \|x^k-x^{k,*}\|
    &\le\frac{\|\zeta\|}{\rho}\le
    \frac{
        \|\nabla f(x^k)\|
        +\beta_k\|J_g(x^k)\|
        +\|\xi^k\|
    }{\rho} \le
    \frac{l_f+\bar\beta l_g+l_r}{\rho},
\end{align*} where the last inequality follows from~\cref{lem:barbeta}
and the assumed bounds on $\nabla f$, $J_g$, and $\partial r$. Combining this  with~\eqref{eq:bndxk-xk*},
using $\rho\ge\underline\rho$, and $1/\ln(1+u)\le \max\left\{1/\ln 2,2/u\right\}$ for all $u>0$, implies that \cref{alg:acceleratedNesterov} will terminate within $T_1^k$ iterations.
Finally, the optimality of the proximal subproblem gives
$
    -\nabla\widehat f(y^t)
    -M_t(x^{t+1}-y^t)
    \in\partial\widehat\Psi(x^{t+1}),
$
and hence
$
    \nabla\widehat f(x^{t+1})
    -\nabla\widehat f(y^t)
    -M_t(x^{t+1}-y^t)
    \in\partial\widehat\phi(x^{t+1});
$
see~\cite[Eq.~(4.10)]{nesterov2013gradient}.
Thus, upon termination, the returned point satisfies
$
    \operatorname{dist}
    \bigl(0,\partial F^k(x^{t+1};\rho)\bigr)
    \le\epsilon.
$
\end{proof}

\section{Cutting-plane methods for \cref{alg:pdmcp}}\label{app:cps}
 This section describes how to use   Vaidya's volumetric cutting-plane method or the ellipsoid method for~\cref{alg:pdmcp}. 

\subsection{Vaidya's method}
\label{app:vaidya}
 We describe how Vaidya's volumetric cutting-plane method~\cite{anstreicher1997vaidya} updates the localization polytope $\mathcal P^t$ and computes the next query point~$y^{t+1}$ in~\cref{alg:pdmcp}. The separation vector $w^t$ is supplied by~\cref{alg:pdmcp} and is taken as input. 
Let
\(
    \mathcal P=\{y\in\mathbb{R}^n:Ay\ge b\}
\)
be a polytope with nonempty interior, where
$A\in\mathbb{R}^{m\times n}$ and $b\in\mathbb{R}^m$.
For $y\in\operatorname{int}(\mathcal P)$, define 
\begin{align}\label{eq:def-sS}
    s(y):=Ay-b,\qquad
    S(y):=\operatorname{Diag}(s(y)),\qquad
    M(y):=A^\top S(y)^{-2}A,
\end{align}
and the volumetric barrier
$
    V_\mathcal{P}(y):=\frac12\ln\det M(y).
$
Let
\[
    \Pi(y):=
    S(y)^{-1}A M(y)^{-1}A^\top S(y)^{-1},
    \qquad
    p_i(y):=\Pi_{ii}(y),\, \forall\, i,
\]
and define
$
    p_{\min}(y):=\min_i p_i(y),$ $
    \Sigma(y):=\operatorname{Diag}(p_1(y),\ldots,p_m(y)),
$ and
\[
    Q(y):=
    A^\top S(y)^{-1}\Sigma(y)S(y)^{-1}A.
\]
For a  symmetric positive definite matrix $Q$, define
$\|z\|_Q:=\sqrt{z^\top Qz}$.

We initialize the procedure with $\mathcal P^0=[-1,1]^n$ and its volumetric center $y^0=0$.
At each nonterminating iteration of~\cref{alg:pdmcp}, \cref{alg:vaidya-update} updates $(\mathcal P^t,y^t)$
to $(\mathcal P^{t+1},y^{t+1})$. Fix constants
$p_{\rm drop}\in(0,1)$, $\tau>0$, and $0<c_1\le c_2$ satisfying the admissibility conditions
in~\cite[Section~6]{anstreicher1997vaidya}, and define the recentering residual
\[
    \mathcal{R}_\mathcal{P}(y):=
    \left\|
        Q(y)^{-1}\nabla V_\mathcal{P}(y)
    \right\|_{Q(y)}
    -
    \min\left\{
        c_1,\,
        \bigl(2\sqrt{p_{\min}(y)}-p_{\min}(y)\bigr)^{1/2}c_2
    \right\}.
\]
We call $y$ an approximate volumetric center of $\mathcal P$
if $\mathcal R_{\mathcal P}(y)\le 0$.
Under the above   conditions, the recentering step
in~\cref{alg:vaidya-update} reaches an approximate volumetric
center in a constant number of Newton steps;
see~\cite[Corollary~6.6]{anstreicher1997vaidya}.

\begin{algorithm}[t!]
\caption{Vaidya Localization Update for \cref{alg:pdmcp}}
\label{alg:vaidya-update} 
\begin{algorithmic}[1]
\REQUIRE
$\mathcal P^t=\{y:A^t y\ge b^t\}$, an interior point $y^t$ satisfying $\mathcal R_{\mathcal P^t}(y^t)\le 0$, 
and a nonzero   vector $w^t$ generated
by~\cref{alg:pdmcp}.
\STATE Using $A=A^t$ and $b=b^t$, compute
$M^t:=M(y^t)$, $p_i^t:=p_i(y^t)$ for all $i$,
and $p_{\min}^t:=\min_i p_i^t$. 

\IF{$p_{\min}^t\ge p_{\rm drop}$}
    \STATE Set $\widetilde b^t := (w^t)^\top y^t -  \sqrt{(w^t)^\top(M^t)^{-1}w^t/\tau}$, 
     $A^{t+1} = [ A^t; (w^t)^\top ]$, and $b^{t+1} = [b^t; \widetilde b^t]$. 
\ELSE
   \STATE Choose $j_t\in\arg\min_i p_i^t$ and form
    $(A^{t+1},b^{t+1})$ by deleting the $j_t$th inequality
    from $A^t y\ge b^t$. 
\ENDIF

\STATE Set  $\mathcal{P}^{t+1}:=\{y:A^{t+1}y\ge b^{t+1}\}$. Starting from $y^t$, apply Newton steps to
$V_{\mathcal P^{t+1}}$ until
$\mathcal R_{\mathcal P^{t+1}}(y)\le 0$.  
Set $y^{t+1}$ to the resulting point.
\RETURN
$(\mathcal{P}^{t+1},y^{t+1})$.
\end{algorithmic}
\end{algorithm}

The update in  \cref{alg:vaidya-update} satisfies
\begin{align}\label{inclusion}
    \mathcal P^t\cap
    \{y:\langle w^t,y-y^t\rangle\ge0\}
    \subseteq \mathcal P^{t+1}.
\end{align}
Indeed, when a new constraint is added,
\( \widetilde b^t < (w^t)^\top y^t, \)
so
$
    \langle w^t,y-y^t\rangle\ge0
$ yields $
    (w^t)^\top y\ge (w^t)^\top y^t>\widetilde b^t.
$
Thus, the halfspace defined by the central cut is contained
in the halfspace defined by the added inequality.
When a constraint is deleted, we have
$\mathcal P^t\subseteq\mathcal P^{t+1}$,
so~\cref{inclusion}   holds as well.

\subsection{Ellipsoid method}\label{app:ellip}
When the ellipsoid method \cite{bland1981ellipsoid} is used for \cref{alg:pdmcp},  the localization set is an ellipsoid whose center serves as the query point. 
For $n\ge 2$, write
$
    \mathcal P^t=\{y^t+A^tu:\|u\|\le 1\},
$
where $A^t\in\RR^{n\times n}$ is nonsingular.
We initialize $y^0=0$ and $A^0=I$, so
$\mathcal P^0=\mathcal B_1$ contains the entire dual domain.
At each nonterminating iteration of~\cref{alg:pdmcp}, the separation vector $w^t$ defines a central cut whose
retained halfspace contains $\mathcal B_{\eta_k}(\bar y)$
by~\cref{lem:ctb}.
The minimum-volume ellipsoid containing the retained
half-ellipsoid is given by the closed-form update
in~\cref{alg:ellipsoid-dcp}.

For $n=1$, $\mathcal P^t$  is an interval
with center $y^t$ and radius $A^t>0$.
In this case, the update reduces to interval bisection:
$
    y^{t+1}=y^t+({A^t}/{2})\operatorname{sign}(w^t),$ and $
    A^{t+1}= {A^t}/{2}.
$
Thus, each update halves the interval length. The update in \cref{alg:ellipsoid-dcp} also guarantees \cref{inclusion}, since $\mathcal P^{t+1}$ contains 
$
    \mathcal P^t\cap
    \{y:\langle w^t,y-y^t\rangle\ge 0\}.
$

\begin{algorithm}[!htbp]
\caption{Ellipsoid Localization Update for \cref{alg:pdmcp}}
\label{alg:ellipsoid-dcp}
\begin{algorithmic}[0]
\REQUIRE
$\mathcal{P}^t=\{y^t+A^tu:\|u\|\le1\}$, its center $y^t$, and a nonzero vector $w^t$ generated by \cref{alg:pdmcp}.
\IF{$n=1$}
    \STATE Set
    $y^{t+1}:=y^t+\frac{A^t}{2}\operatorname{sign}(w^t)$. Set $A^{t+1}:=A^t/2$.
\ELSE 
\STATE Set
    $q^t:=(A^t)^\top w^t/\|(A^t)^\top w^t\|$. Set
    $y^{t+1}:=y^t+A^tq^t/(n+1)$.
    \STATE Set $A^{t+1}=\sqrt{n^2/(n^2-1)}(A^t-\bigl(1-\sqrt{(n-1)/(n+1)}\bigr)A^tq^t(q^t)^\top)$.
\ENDIF
\STATE Set $\mathcal P^{t+1}=\{y^{t+1}+A^{t+1}u:\|u\|\le1\}$.
\RETURN
$(\mathcal{P}^{t+1},y^{t+1})$.
\end{algorithmic}
\end{algorithm}

\FloatBarrier
\bibliographystyle{siamplain}
\bibliography{optim}

\end{document}